\documentclass[11pt,twoside]{article}  

\usepackage{fullpage, amsthm, amsmath, amssymb, hyperref, dsfont, mathtools,bm}
\usepackage{framed}
\usepackage{caption,subcaption,graphicx,rotating,multirow}
\usepackage{hyperref,url}
\usepackage{color,xcolor}
\usepackage{enumitem}
\usepackage{epstopdf}
\usepackage[numbers]{natbib}

\theoremstyle{plain}
\newtheorem{theorem}{Theorem}

\newtheorem{lemma}[theorem]{Lemma}

\newtheorem{corollary}[theorem]{Corollary}

\newtheorem{proposition}[theorem]{Proposition}
 \newenvironment{proofof}[1]{{\bf {\em Proof of #1.}}}{\hfill \rule{2mm}{2mm} 
 }
  
\newlength{\widebarargwidth}
\newlength{\widebarargheight}
\newlength{\widebarargdepth}
\DeclareRobustCommand{\widebar}[1]{%
  \settowidth{\widebarargwidth}{\ensuremath{#1}}%
  \settoheight{\widebarargheight}{\ensuremath{#1}}%
  \settodepth{\widebarargdepth}{\ensuremath{#1}}%
  \addtolength{\widebarargwidth}{-0.3\widebarargheight}%
  \addtolength{\widebarargwidth}{-0.3\widebarargdepth}%
  \makebox[0pt][l]{\hspace{0.3\widebarargheight}%
    \hspace{0.3\widebarargdepth}%
    \addtolength{\widebarargheight}{0.3ex}%
    \rule[\widebarargheight]{0.95\widebarargwidth}{0.1ex}}%
  {#1}}
  
\makeatletter
\long\def\@makecaption#1#2{
        \vskip 0.8ex
        \setbox\@tempboxa\hbox{\small {\bf #1:} #2}
        \parindent 1.5em  
        \dimen0=\hsize
        \advance\dimen0 by -3em
        \ifdim \wd\@tempboxa >\dimen0
                \hbox to \hsize{
                        \parindent 0em
                        \hfil 
                        \parbox{\dimen0}{\def\baselinestretch{0.96}\small
                                {\bf #1.} #2
                                } 
                        \hfil}
        \else \hbox to \hsize{\hfil \box\@tempboxa \hfil}
        \fi
        }
\makeatother

\newcommand{\kull}[2]{\ensuremath{D(#1\; \| \; #2)}}

\newcommand{\unif}{\ensuremath{\mbox{Unif}}}

\newcommand{\Fspace}{\ensuremath{\mathcal{F}}}

\newcommand{\bigo}{\ensuremath{\mathcal{O}}}

\newcommand{\E}{\operatorname{\mathbb{E}}}

\DeclareMathOperator{\vol}{vol}

\newcommand{\Z}{\mathbb{Z}}
\newcommand{\Zp}{\mathbb{Z}_{+}}
\newcommand{\R}{\mathbb{R}}

\newcommand{\Rp}{\mathbb{R}_{+}}

\newcommand{\indi}{\mathds{1}}

\newcommand{\imnb}{\mathbf{i}}

\newcommand{\twonorm}[1]{\left\|#1\right\|_{\ell_2}}

\newcommand{\abs}[1]{\left|#1\right|}
\newcommand{\vct}[1]{\bm{#1}}
\newcommand{\mtx}[1]{\underline{\bm{#1}}}

\def\BC{\begin{center}}
\def\EC{\end{center}}
\def\BIT{\begin{itemize}}
\def\EIT{\end{itemize}}
\def\BET{\begin{enumerate}}
\def\EET{\end{enumerate}}
\def\BEQ{\begin{equation}}
\def\EEQ{\end{equation}}

\long\def\comment#1{}

\usepackage{tikz}
\usepackage{tkz-graph}					
\usepackage[margin = 2.5cm]{geometry}

\newcommand{\numvar}{\ensuremath{n}}
\newcommand{\numtype}{\ensuremath{d}}
\newcommand{\numobs}{\ensuremath{m}}

\newcommand{\one}{\mathbf{1}}
\newcommand{\half}{\frac{1}{2}}

\newcommand{\Zpart}{\mathcal{Z}}
\newcommand{\fenergy}{\mathfrak{F}}
\newcommand{\overlap}{\mu}
\newcommand{\matoverlap}{\mtx{\overlap}}
\newcommand{\eulerflow}{\nu}
\newcommand{\mateulerflow}{\mtx{\eulerflow}}
\newcommand{\ncc}{\operatorname{\mathsf{ncc}}} 
\newcommand{\nst}{\operatorname{\mathsf{nst}}} 
\newcommand{\erfc}{\operatorname{\mathsf{erfc}}} 

\newcommand{\BP}{\mathsf{BP}} 
\newcommand{\IT}{\mathsf{IT}} 
\newcommand{\HQP}{\mathsf{HQP}} 

\newcommand{\x}{\underline{\bm{x}}}
\newcommand{\y}{\underline{\bm{y}}} 
\newcommand{\z}{\underline{\bm{z}}}
\newcommand{\w}{\underline{\bm{w}}}

\newcommand{\probquery}{\alpha}
\newcommand{\typeprop}{\pi}

\usepackage{etoolbox}
\makeatletter
\patchcmd{\maketitle}{\@fnsymbol}{\@alph}{}{}  
\makeatother

\title{
{\bf{\LARGE{Decoding from Pooled Data: \\ Sharp Information-Theoretic Bounds}}}
}
\author{Ahmed El Alaoui\thanks{Department of Electrical Engineering and Computer Sciences, UC Berkeley, CA.} ~~ 
Aaditya Ramdas$^{\star}$\thanks{Department of Statistics, UC Berkeley, CA.} \and
Florent Krzakala\thanks{Laboratoire de Physique Statistique, CNRS, PSL Universit\'es \& Ecole Normale Sup\'erieure, Sorbonne Universit\'es et Universit\'e Pierre \& Marie Curie, Paris, France.} ~~ 
Lenka Zdeborov\'{a}\thanks{Institut de Physique Th\'eorique, CNRS, CEA, Universit\'e Paris-Saclay, Gif-sur-Yvette, France.} ~~ 
Michael I. Jordan$^{\star \dagger}$}

\begin{document}
\date{}
\maketitle


\vspace*{-.3in} 

\begin{abstract}

Consider a population consisting of $n$ individuals, each of whom has one of $d$ types (e.g. their blood type, in which case $d=4$). We are allowed to query this database by specifying a subset of the population, and in response we observe a noiseless histogram (a $d$-dimensional vector of counts) of types of the pooled individuals. This measurement model arises in practical situations such as pooling of genetic data and may also be motivated by privacy considerations. We are interested in the number of queries one needs to unambiguously determine the type of each individual. 
In this paper, we study this information-theoretic question under the random, dense setting where in each query, a random subset of individuals of size proportional to $n$ is chosen. This makes the problem a particular example of a random constraint satisfaction problem (CSP) with a ``planted'' solution. We establish almost matching upper and lower bounds on the minimum number of queries $m$ such that there is no solution other than the planted one with probability tending to 1 as $n \to \infty$. 
Our proof relies on the computation of the exact ``annealed free energy" of this model in the thermodynamic limit, which corresponds to the exponential rate of decay of the expected number of solution to this planted CSP. As a by-product of the analysis, we show an identity of independent interest relating the Gaussian integral over the space of Eulerian flows of a graph to its spanning tree polynomial.
\end{abstract}

\section{Introduction}
Constraint satisfaction problems (CSPs) have been the object of intense study
in recent years in probability theory, computer science, information theory and
statistical physics.  For certain families of CSPs, a deep understanding has begun
to emerge regarding the number of solutions as a function of problem size, as well as the
algorithmic feasibility of finding solutions when they exist (see e.g.
\cite{coja2009spectral,coja2014analyzing,coja2016walksat,coja2016belief,ding2015proof,ding2016satisfiability,sly2016number}...)
Consider in particular a \emph{planted} random constraint satisfaction problem
with $\numvar$ variables that take their values in the discrete set $\{1,\cdots,\numtype\}$,
with $\numtype \ge 2$, and with $\numobs$ clauses drawn uniformly at random under the
constraint that they are all satisfied by a pre-specified assignment, which is referred
to as \emph{the planted solution}. It is of interest to determine the properties
of the set of all solutions of CSP as $\numvar$ and $\numobs$ grow to infinity at a some
relative rate. Two questions are of particular importance: \textit{(1) how large should
$\numobs$ be so that the planted solution is the unique solution to CSP?} and
\textit{(2) given that it is unique, how large should $\numobs$ be so that it is
recoverable by a ``tractable" algorithm?} Significant progress has been made on
these questions, often initiated by insights from statistical physics and followed
by a growing body of rigorous mathematical investigation. The emerging picture is
that in many planted CSPs, when $\numvar$ is sufficiently large, all solutions
become highly correlated with the planted one when $\numobs > \kappa_{\IT}\cdot \numvar$,
for some ``Information-Theoretic" ($\IT$) constant $\kappa_{\IT} > 0$. Furthermore,
one of these highly correlated solutions becomes typically recoverable 
by a random walk or a Belief Propagation ($\BP$)-inspired algorithm when 
$\numobs > \kappa_{\BP}\cdot \numvar$ for some $\kappa_{\BP} >
\kappa_{\IT}$~\cite{coja2009spectral,krzakala2009hiding,krzakala2012reweighted,coja2014analyzing}.
Interestingly, it is known in many problems, at least heuristically, that these algorithms 
fail when $\kappa_{\IT} < \numobs/\numvar < \kappa_{\BP}$, and a tractable algorithm
that succeeds in this regime is still lacking~\cite{achlioptas2008algorithmic,coja2009random,zdeborova2015statistical,coja2016walksat}.
In other words, there is a non-trivial regime $\numobs/\numvar \in (\kappa_{\IT},\kappa_{\BP})$
where an essentially unique solution exists, but is hard to recover. In the picture we
described, the uniqueness and recoverability thresholds differ by a constant factor;
that is, they both happen at the same scale where $\numobs$ is proportional to $\numvar$.
In other examples, notably a class of planted CSPs related to XORSAT, the gap
between the information theoretic and the putative algorithmic threshold is much
larger~\cite{feldman2015complexity}.

The systematic study of this phenomenon has mostly dealt with planted CSPs with
Boolean variables (the $\numvar$ variables are $\{0,1\}$-valued and $\numtype =2$).
Among CSPs with larger variable domain, the random Graph Coloring problem with
$\numtype$ colors is a prominent prototype.  Here the size of the gap is constant.


In this paper we consider a naturally motivated random CSP with arbitrary but fixed domain size $\numtype$, which we call the Histogram Query Problem ($\HQP$), where we present evidence that the gap between existence of a unique solution and its tractable recovery could be as large as $\log n$. More precisely, in this paper, we undertake a detailed information-theoretic analysis which shows that in $\HQP$, the planted solution becomes unique at as soon as $\numobs > \gamma^* \numvar/\log \numvar$ with high probability as $\numvar \to \infty$ for a constant $\gamma^* > 0$. 
In a sequel paper, we consider the algorithmic aspect of the problem and provide a $\BP$-based algorithm that  recovers the planted assignment if $\numobs \ge \kappa^*\cdot \numvar$ for a specific threshold $\kappa^*$ and fails otherwise. This indicates the existence of a logarithmic gap between the information-theoretic and algorithmic thresholds. 

\subsection{Problem and motivation}

\paragraph{The setting}

Let $\{\vct{h}_a\}_{1\le a \le m}$ be a collection of $\numtype$-dimensional arrays with non-negative integer entries. For an assignment $\tau  : \{1,\cdots,\numvar\} \mapsto \{1,\cdots,\numtype\}$ of the $\numvar$ variables, and given a realization of $\numobs$ random subsets $S_a \subset \{1,\cdots,\numvar\}$, the constraints of the $\HQP$ are given by $\vct{h}_a = \vct{h}_a(\tau)$ for all $1\le a \le \numobs$ with 
\[\vct{h}_a(\tau) := \left(\left|{\tau}^{-1}(1) \cap S_a \right|,\cdots, \left|{\tau}^{-1}(\numtype) \cap S_a \right| \right) \in \Zp^{\numtype}.\] 
We let $\tau^* : \{1,\cdots,\numvar\} \mapsto \{1,\cdots,\numtype\}$ be a planted assignment; i.e., we set $\vct{h}_a := \vct{h}_a(\tau^*)$ for all $a$ for some realization of the sets $\{S_a\}$, and consider the problem of recovering the map $\tau^*$ given the observation of the arrays $\{\vct{h}_{a}\}_{1\le a \le \numobs}$. 

This problem can be viewed informally as that of decoding a discrete high-dimensional signal consisting of categorical variables from a set of measurements formed by pooling together the variables belonging to a subset of the signal. It is useful to think of the $\numvar$ variables as each describing the type or category of an individual in a population of size $\numvar$, where each individual has exactly one type among $\numtype$. For instance the categories may represent blood types or some other discrete feature such as ethnicity or age group. Then, the observation $\vct{h}_a$ is the histogram of types of a subpopulation $S_a$. We let $\vct{\pi} = \frac{1}{n} \left(\left|{\tau^*}^{-1}(1)\right|,\cdots,\left|{\tau^*}^{-1}(d)\right| \right)$ denote the vector of proportions of assigned values; i.e., the empirical distribution of categories.

We consider here a model in which each variable participates in a given constraint independently and with probability $\alpha\in (0,1)$. Thus, the sets $\{S_a\}_{1\le a \le \numobs}$ are independent draws of a random set $S$ where $\Pr(i \in S) = \alpha$ independently for each $i \in \{1,\ldots,n\}$. We are thus in the ``dense regime" where $\E[|S|] = \alpha \numvar$; i.e., the number of variables participating in each constraint (the degree of each factor in the CSP) is linear in $\numvar$. 

\paragraph{Motivation}
This model is inspired by practical problems in which a data analyst can only assay certain summary statistics involving a moderate or large number of participants. This may be done for privacy reasons, or 
it may be inherent in the data-collection process (see e.g.\ \cite{sham2002dna,Heo1163}).  For example, in DNA assays, the pooling of allele measurements across multiple strands of DNA is necessary given the impracticality of separately analyzing
individual strands.  Thus the data consists of a frequency spectrum of alleles; a ``histogram'' in our language.  In the privacy-related situation, one may take the viewpoint of an attacker whose goal is to gain a granular knowledge of the database from coarse measurements, or that of a guard who wishes to prevent this scenario from happening. It is then natural to ask how many histogram queries it takes to exactly determine the category of each individual.

\paragraph{Related problems}
Note that the case $\numtype=2$ of $\HQP$ can be seen as a compressed sensing problem with a binary sensing matrix and binary signal. While the bulk of the literature in the field of compressed sensing is devoted to the case in which both the signal of interest and the sensing matrix are real-valued, the binary case has also been considered, notably in relation to Code Division Multiple Access (CDMA) \cite{zigangirov2004theory,tanaka2002statistical}, and Group Testing \cite{du2006pooling,mezard2011group}: in the latter, one observes the logical ``OR" of subsets of the entries of the signal. In the case of categorical variables with $\numtype \ge 3$ categories, it is natural to consider measurements consisting of histograms of the categories in the pooled sub-population. 
In the literature on compressed sensing one commonly considers the setting where the sensing matrices have i.i.d.\ entries with finite second moment, and the signal has an arbitrary empirical distribution of its entries. It has been established that, under the scaling $\numobs = \kappa \numvar$, whereas the success of message-passing algorithms requires $\kappa>\kappa_{\BP}$ \cite{bayati2015universality}, the information-theoretic threshold is $\kappa_{\IT}=0$ in the discrete signal case~\cite{wu2009fundamental,donoho2013information}, indicating that uniqueness of the solution happens at a finer scale $\numobs = o(\numvar)$. Here we consider the $\HQP$ with arbitrary $\numtype$, for which the exact scaling for investigating uniqueness is $\numobs = \gamma \frac{\numvar}{\log \numvar}$ with finite $\gamma >0$, and provide tight bounds on the information-theoretic threshold.     

\paragraph{Prior work on $\HQP$}
The study of this problem for generic values of $d$ was initiated in~\cite{wang2016data} in the two settings where the sets $\{S_a\}$ are deterministic and random. They showed in both these cases with a simple counting argument that under the condition that $\vct{\pi}$ is the uniform distribution, if $m < \frac{\log d}{d-1}\frac{n}{\log n}$ then the set of collected histograms does not uniquely determine the planted assignment $\tau^*$ (with high probability in the random case). On the other hand, for the deterministic setting, they provided a querying strategy that recovers $\tau^*$ provided that $m > c_0 \frac{n}{\log n}$, where $c_0$ is an absolute constant independent of $d$.  For the random setting and under the condition that the sets $S_a$ are of average size $\numvar/2$, they proved via a first moment bound that $m > c_1 \frac{n}{\log n}$ with $c_1$ also constant and independent of $d$, suffices to uniquely determine $\tau^*$, although no algorithm was proposed in this setting. 

\vspace{.2cm}
In the above results, there is a gap that is both information-theoretic and algorithmic depending on the dimension $\numtype$ between the upper and lower bounds. Intuitively, the upper bounds should also depend on $\numtype$ since the decoding problem becomes easier (or at least, it is no harder) for large $d$, for the simple reason that if it is possible to determine the categories of the population for $\numtype=2$, then one can proceed by dichotomy for larger $\numtype$ by merging the $\numtype$ groups into two super-groups, identifying which individuals belong to each of the two super-groups, and then recurse. We attempt to fill the information-theoretic gap in the random setting by providing tighter upper and lower bounds on the number of queries $\numobs$ necessary and sufficient to uniquely determine the planted assignment $\tau^*$ with high probability, which depend on the dimension $\numtype$ and $\vct{\pi}$ along with explicit constants. In a sequel paper, we consider the algorithmic aspect of the problem and provide a Belief Propagation-based algorithm that recovers the planted assignment if $\numobs \ge \kappa^*(\vct{\pi},d)\cdot \numvar$ for a specific threshold $\kappa^*(\vct{\pi},d)$ and fails otherwise, indicating the putative existence of a statistical-computational gap in the random setting.       

\subsection{Main result} 
Let $\Delta^{d-1}$ be the $d-1$-dimensional simplex and $H(\vct{x}) = - \sum_{r=1}^d x_{r}\log x_{r}$ for $\vct{x} \in \Delta^{d-1}$ be the Shannon entropy function. We write $\tau \sim \vct{\pi}$ to indicate that $\tau$ is a random assignment drawn from the uniform distribution  over maps $\tau: \{1,\cdots,n\} \mapsto \{1,\cdots,d\}$ such that ${\frac{1}{n} \left(\left|{\tau}^{-1}(1)\right|,\cdots,\left|{\tau}^{-1}(d)\right| \right) = \vct{\pi}}$. 
\begin{theorem} \label{mainthm}
For $n\ge 2$ integer, $m = \gamma \frac{n}{\log n}$, $\gamma >0$, $\alpha \in (0,1)$, and $\vct{\pi} \in \Delta^{d-1}$ with entries bounded away from 0 and 1.
Let $\mathcal{E}$ be the event that $\tau^*$ is \emph{not} the unique satisfying assignment to $\HQP$:
\[\mathcal{E} = \left\{\exists \tau \in \{1,\cdots,\numtype\}^\numvar~:~ \tau \neq \tau^*, 
~ 
\vct{h}_a(\tau) = \vct{h}_a(\tau^*)
 ~ \forall a \in \{1,\cdots,\numobs\} 
\right\}.\] 
{\bf (i)} If 
\[\gamma < \gamma_{\text{low}} := \frac{H(\vct{\pi})}{d-1},\]
then 
\[\lim_{n \to \infty} \E_{\tau^* \sim\vct{\pi}}\Pr \left(\mathcal{E} \right) = 1.\]
{\bf (ii)} On the other hand, let $\vct{\pi}_{[\cdot]}$ be the vector of order statistics of $\vct{\pi}$: $\typeprop_{[1]} \ge \pi_{[2]} \ge \cdots \ge \pi_{[d]}$. For $1\le k \le d-1$, let $\vct{\pi}^{(k)} \in \Delta^{k-1}$ be defined as 
${\pi}^{(k)}_{1} = \sum_{r=1}^{d-k+1} \pi_{[r]}$ and ${\pi}^{(k)}_{l} = \pi_{[d-k+l]}$ for all $2\le l \le k$ (if $k\ge 2$). 
If 
\[\gamma > \gamma_{\text{up}} := 2\max_{ 1 \le k \le d-1} \frac{H(\vct{\pi}) - H({\vct{\pi}}^{(k)})}{d-k},\] 
then 
\[\lim_{n \to \infty} \E_{\tau^* \sim\vct{\pi}}\Pr\left(\mathcal{E} \right) = 0.\]

\end{theorem}
 
 \paragraph{Remarks and special cases:} 
\begin{itemize}
\item For $d=2$, $\gamma_{\text{up}} = 2 H(\vct{\pi}) = 2\gamma_{\text{low}}$.
\item If $\vct{\pi} = (\frac{1}{d},\cdots,\frac{1}{d})$, or more generally, if $\vct{\pi}$ is such that $k=1$ maximizes the expression defining $\gamma_{\text{up}}$ then $\gamma_{\text{up}} = 2 \frac{H(\vct{\pi})}{d-1} = 2\gamma_{\text{low}}$.
 \item The resulting bounds do not depend on $\alpha$ as long as it is fixed and bounded away from 0 and 1. Its contribution in the problem is sub-dominant and vanishes as $n \to \infty$ under the scaling considered here.  
\item The number $k$ in the expression of $\gamma_{\text{up}}$ can be interpreted as the number of connected components of a graph on $d$ vertices that depends on the overlap structure of the two assignments $\tau$ and $\tau^*$, and induces ``maximum confusion" between them. This will become clear in latter sections. 
\end{itemize}

The proof of the above Theorem occupies the rest of the manuscript.

\subsection{Main ideas of the proof} 
Our main contribution is the second part of Theorem~\ref{mainthm}, which establishes an upper
bound on the uniqueness threshold of the random CSP with histogram constraints $\HQP$.
The proof uses the first moment method to upper bound the probability of existence
of a non-planted solution. Since we are in a planted model, the analysis of the first
moment ends up bearing many similarities with a second moment computation in a purely
random (non-planted) model. Although second moment computations often require approximations,
for the $\HQP$ it turns out that we are able to compute the exact annealed free energy
of the model in the thermodynamic limit. That is, letting $\Zpart$ be the number of
solutions of the CSP, we show that the limit
\[\fenergy(\gamma) := \lim_{n \to \infty} \frac{1}{n} \log \E\left[\Zpart-1\right]\]
exists and we compute its value exactly.
Then the value of the threshold $\gamma_{\text{up}}$ is obtained by locating the first
point at which $\fenergy$ becomes negative:
\[\gamma_{\text{up}} = \inf\left\{\gamma>0 ~:~ \fenergy(\gamma) <0\right\}.\]
Together with the fact that $\fenergy$ is a monotone function, which will become
clear once $\fenergy$ is computed, it is clear that for any $\gamma > \gamma_{\text{up}}$,
$\E[\Zpart-1]$ decays exponentially with $\numvar$ when the latter is sufficiently large.

This general strategy has been successfully pursued for a range of CSPs, such as K-SAT, NAE-SAT, and Independent Set, most of which are Boolean. 
For larger domain sizes, in order to carry out the second moment method one needs fine control of the overlap structure
between the planted and a candidate solution. 
This control is at the core of the difficulty that arises in any second moment computation.
 To obtain such control, researchers have often imposed additional assumptions, 
 at a cost of a weakening of the resulting bounds. 
For example, existing proofs for Graph Coloring and similar problems assume certain balancedness conditions
(the overlap matrix needs to be close to doubly stochastic.) without which the annealed
free energy cannot be computed~\cite{achlioptas2004two_journal,achlioptas2004chromatic,coja2016chromatic,
bapst2016condensation,banks2016information}; this yields results that fall somewhat short of
the bounds that the second moment method could achieve in principle~\cite{dani2012tight}.
In the present problem, due its rich combinatorial structure, we are able to obtain
unconditional control of the overlap structure, for any domain size $d$, and compute
the exact annealed free energy.

Concretely, computing the function $\fenergy$ requires tight control of
the ``collision probability" of two non-equal assignments $\tau_1$ and $\tau_2$.
This is the probability that the random histograms $\vct{h}(\tau_1) = \left(\left|{\tau_1}^{-1}(1) \cap S \right|,\cdots, \left|{\tau_1}^{-1}(\numtype) \cap S \right| \right)$
and $\vct{h}(\tau_2) = \left(\left|{\tau_2}^{-1}(1) \cap S \right|,\cdots, \left|{\tau_2}^{-1}(\numtype) \cap S \right| \right)$ generated from a random draw of a pool $S$ coincide. The collision probability roughly measures the correlation strength between the two assignments.
Specifically, we will be interested in the collision probabilities of
the pairs $(\tau^*,\tau)$ where $\tau^*$ is the planted assignment and $\tau$ is any
candidate assignment. Its decay reveals how long an assignment $\tau$ ``survives"
as a satisfying assignment to $\HQP$ as $n \to \infty$.  The study of these collision
probabilities requires the evaluation of certain Gaussian integrals over the space
of \emph{Eulerian flows} of a weighted graph on $\numtype$ vertices that is defined
based on the overlap structure of $\tau$ and $\tau^*$. We prove a family of identities
that relate these integrals to some combinatorial polynomials in the weights of the graph: 
the spanning tree and spanning forest
polynomials. We believe that these identities are of independent
interest beyond the problem studied in this paper. Once these collision probabilities
are controlled, the computation of $\fenergy(\gamma)$ per se requires the analysis of
a certain sequence of optimization problems. We show that the sequence of maximum
values converges to a finite limit that yields the value of the annealed free energy.

On the other hand, the proof of the first part of Theorem~\ref{mainthm} is straightforward---it is an extension
of a standard counting argument used in~\cite{zhang2013non} and~\cite{wang2016data}. The argument
goes as follows: if $\numobs$ is too small then the number of possible histograms one
could potentially observe is exponentially smaller than the number of assignments of
$n$ variables that agree with $\vct{\pi}$. Therefore when the planted assignment
$\tau^*$ is drawn at random, there will exist at least one $\tau \neq \tau^*$ that
satisfies the constraints of the CSP with overwhelming probability.  We begin with
this argument in the next section and then turn to the more challenging
computation of the upper bound.

\section{Proof of Theorem~\ref{mainthm}}
\paragraph{Notation} We denote vectors in $\R^d$ in bold lower case letters, e.g., $\vct{x}$, and matrices in $\R^{d \times d}$ will be written in bold lower case underlined letters, e.g., $\mtx{x}$. We denote  the coordinates of such vectors and matrices as $x_r$ and $x_{rs}$ respectively. Matrices that act either as linear operators on the space $\R^{d\times d}$ or that are functions of elements in this space are written in bold upper case letters, e.g., $\bm{M}\x$ and $\bm{L}(\x)$, for $\x \in \R^{d \times d}$.  These choices will be clear from the context. We may write $\x/\y$ to indicate coordinate-wise division. Additionally, for two $d \times d$ matrices $\mtx{a}, \mtx{b} \in \R^{d \times d}$, $\mtx{a} \odot \mtx{b} \in \R^{d \times d}$ is their Hadamard product. We let $\one \in \R^d$ be the all-ones vector. %
\subsection{The first part of Theorem~\ref{mainthm}: the lower bound} 
Let $\numobs = \gamma \frac{\numvar}{\log \numvar}$ with $\gamma >0$. The number of potential histograms one could possibly observe in a single query with pool size $|S| = k$ is $f(k,\numtype) := {{\numtype+k-1}\choose{\numtype-1}} \le (k+1)^{\numtype-1}$. Since the queries are independent, the number of collections of histograms $\{\vct{h}_a\}_{1\le a \le \numobs}$ one could potentially observe in $\numobs$ queries is $\prod_{a=1}^m f(|S_a|,\numtype)$.
On the other hand, the number of possible assignments $\tau : \{1,\cdots,\numvar\} \mapsto \{1,\cdots,\numtype\}$ satisfying the constraint $\vct{\pi} = \frac{1}{n} \left(\left|{\tau^*}^{-1}(1)\right|,\cdots,\left|{\tau^*}^{-1}(d)\right| \right)$ is ${{\numvar}\choose{\numvar\vct{\pi}}}={{\numvar}\choose{\numvar\pi_1,\cdots,\numvar\pi_\numtype}} \ge C(\vct{\pi})\numvar^{-(\numtype-1)/2} \exp (H(\vct{\pi}) \numvar)$, for some constant $C(\vct{\pi})>0$ depending on $\vct{\pi}$. 

Now, the probability that $\tau^*$ is the unique satisfying assignment of the CSP with constraints given by the random histograms $\{\vct{h}_a(\tau^*)\}_{1\le a \le \numobs}$, averaged over the random choice of $\tau^* \sim \vct{\pi}$, is 
\begin{align*}
\E_{\tau^*\sim \vct{\pi}} \E_{\{S_a\}} &\Big[ \indi\{\forall \tau \in \{1,\cdots,\numtype\}^\numvar ~:~ \vct{h}_a(\tau) = \vct{h}_a(\tau^*) ~\forall a \in \{1, \cdots, \numobs\} ~\implies ~\tau = \tau^*\}\Big] \\
&\le {{\numvar}\choose{\numvar\vct{\pi}}}^{-1} \cdot \E_{S} \left[f(|S|,\numtype)\right]^\numobs\\ 
&\le  {{\numvar}\choose{\numvar\vct{\pi}}}^{-1} \cdot \E_{S}\left[(|S|+1)^{\numtype-1}\right]^\numobs \\
&\le C(\vct{\pi})~\numvar^{(d-1)/2} \cdot \exp\Big(-H(\vct{\pi})\numvar\Big) \cdot (\numvar+1)^{\numobs(\numtype-1)}\\
&\le C(\vct{\pi})~\numvar^{(d-1)/2} \cdot \exp \Big((\gamma(\numtype-1) - H(\vct{\pi}))\numvar\Big).
\end{align*}
If $\gamma < \gamma_{\text{low}}$ the last quantity tends to 0 as $n \to \infty$. This concludes the proof of the first assertion of the theorem.

\subsection{The second part of Theorem~\ref{mainthm} : the upper bound} 
We use a first moment method to show that when $\gamma$ is greater than $\gamma_{\text{up}}$, the only assignment satisfying $\HQP$ is $\tau^*$ with high probability. Let $\Zpart$ be the number of satisfying assignments to $\HQP$:
\begin{equation}\label{Zpartition}
\Zpart := \big|\{\tau \in \{1,\cdots,\numtype\}^\numvar ~:~ \vct{h}_a(\tau) = \vct{h}_a(\tau^*) ~ \forall a \in \{1,\cdots,\numobs\}\}\big|.
\end{equation}
The planted assignment $\tau^*$ is obviously a solution, so we always have $\Zpart \ge 1$.  
Recall the definition of the annealed free energy  
\begin{equation} \label{free_energy_def}
\fenergy(\gamma) := \lim_{n \to \infty} \frac{1}{n} \log \E\left[\Zpart-1\right].
\end{equation}
Also, recall that for $1\le k \le d-1$, $\vct{\pi}^{(k)} \in \Delta^{k-1}$ be defined as 
${\pi}^{(k)}_{1} = \sum_{r=1}^{d-k+1} \pi_{[r]}$ and ${\pi}^{(k)}_{l} = \pi_{[d-k+l]}$ for all $2\le l \le k$ (if $k\ge 2$). 

\begin{theorem}\label{free_energy_value}
Let $m = \gamma \frac{n}{\log n}$ with $\gamma >0$. The limit~\eqref{free_energy_def} exists for all $\gamma >0$ and its value is
\begin{equation}\label{free_value}
\fenergy (\gamma) = \max_{ 1 \le k \le d-1} \left\{ H(\vct{\pi}) - H(\vct{\pi}^{(k)}) - \frac{\gamma}{2}(d-k)\right\}.
\end{equation}
\end{theorem}
We can deduce from Theorem~\ref{free_energy_value} the smallest value of $\gamma$ past which $\fenergy(\gamma)$ becomes negative. In particular, we see that $\fenergy$ is a decreasing function of $\gamma$ that crosses the horizontal axis at 
\[\gamma_{\text{up}} = 2\max_{ 1 \le k \le d-1} \frac{H(\vct{\pi}) - H(\vct{\pi}^{(k)})}{d-k}.\]
From this result it is easy to prove the second assertion of Theorem~\ref{mainthm}. By averaging over $\tau^*$ and applying Markov's inequality, we have:
\[\E_{\tau^* \sim\vct{\pi}}\Pr\left(\exists \tau \in \{1,\cdots,\numtype\}^\numvar~:~ \tau \neq \tau^*, \vct{h}_a(\tau) = \vct{h}_a(\tau^*) ~ \forall a \in \{1,\cdots,\numobs\} \right)  =\E_{\tau^* \sim\vct{\pi}} \Pr\left(\Zpart \ge 2\right) \le \E[\Zpart-1].\] 
For $\gamma > \gamma_{\text{up}}$, it is clear that $\fenergy(\gamma) <0$. Let $0 < \epsilon < \abs{\fenergy(\gamma)}/2$; then there is an integer $n_0(\epsilon)\ge 0$ such that for all $n \ge n_0(\epsilon)$,  
\begin{align*}
\E_{\tau^* \sim\vct{\pi}}\Pr\left(\exists \tau \in \{1,\cdots,\numtype\}^\numvar~:~ \tau \neq \tau^*, \vct{h}_a(\tau) = \vct{h}_a(\tau^*) ~ \forall a \in \{1,\cdots,\numobs\} \right) &\le \exp n ~(\fenergy(\gamma)+\epsilon),\\
&\le \exp n ~\fenergy(\gamma)/2,\\
&\underset{n \to \infty}{\xrightarrow{\qquad}} 0.
\end{align*}
Now it remains to prove Theorem~\ref{free_energy_value}, and this represents the main technical thrust of our paper.

\subsection{Collisions, overlaps, and the first moment}

\paragraph{Preliminaries} We begin by presenting the main quantities to be analyzed in our application of the first moment method. 
We have 
\begin{align*}
\E_{\tau^* \sim \vct{\pi}}\E_{\{S_a\}}[\Zpart - 1] &= \E_{\tau^* \sim \vct{\pi}}\left[\sum_{\substack{\tau \in \{1,\cdots,d\}^n\\ \tau \neq \tau^*}} \Pr\left(\vct{h}_a(\tau) = \vct{h}_a(\tau^*) ~\forall a \in \{1,\cdots,\numobs\} \right) \right] \\
& =  (\numtype^\numvar - 1) \Pr_{\tau,\tau^*,\{S_a\}}\left(\vct{h}_a(\tau) = \vct{h}_a(\tau^*)~\forall a \in \{1,\cdots,\numobs\} \right),
\end{align*}
where $\tau^* \sim \vct{\pi}$, $\tau \sim \unif(\{1,\cdots,\numtype\}^n\backslash\{\tau^*\})$. By conditional independence, 
\begin{align*}
\Pr_{\tau,\tau^*,\{S_a\}}\left(\vct{h}_a(\tau) = \vct{h}_a(\tau^*)~\forall a \in \{1,\cdots,\numobs\} \right) &= \E_{\tau,\tau^*}\left[\Pr_{\{S_a\}}\left(\vct{h}_a(\tau) = \vct{h}_a(\tau^*)~\forall a \in \{1,\cdots,\numobs\} \right)\right] \\
&= \E_{\tau,\tau^*}\left[\Pr_{S}\left(\vct{h}(\tau) = \vct{h}(\tau^*)\right)^m\right].
\end{align*}

Next, we write the \emph{collision probability}, $\Pr_{S}\left(\vct{h}(\tau) = \vct{h}(\tau^*)\right)$, for fixed $\tau$ and $\tau^*$ in a convenient form. Let us first define the \emph{overlap matrix}, $\matoverlap(\tau,\tau^*) = (\overlap_{rs})_{1 \le r,s\le d} \in \Zp^{\numtype \times \numtype}$, of $\tau$ and $\tau^*$, by 
\begin{equation}\label{overlapmat}
\overlap_{rs} = \left|\tau^{-1}(r) \cap {\tau^*}^{-1}(s)\right| \quad \mbox{for all} ~ r,s=1,\cdots,d.
\end{equation}
Remark that $\vct{h}(\tau) = \vct{h}(\tau^*)$ if and only if $\left|S \cap \tau^{-1}(r)\right| = \left|S \cap {\tau^*}^{-1}(r)\right|$ for all $r \in \{1,\cdots,\numtype\}$. Since the collection of sets $\{\tau^{-1}(r)\}_{1\le r \le \numtype}$ forms a partition of $\{1,\cdots,\numvar\}$, and similarly with $\tau^*$, we have the following equality of events 
\[\left\{\vct{h}(\tau) = \vct{h}(\tau^*)\right\} = \left\{ \sum_{s=1}^d \left|S \cap \tau^{-1}(r) \cap {\tau^*}^{-1}(s)\right| = \sum_{s=1}^d \left|S \cap \tau^{-1}(s) \cap {\tau^*}^{-1}(r)\right|,~ \forall r \in \{1,\cdots, \numtype\} \right\}. \]
Therefore, the probability that two assignments $\tau$ and $\tau^*$ collide on a random pool $S$---meaning that their histograms formed on the pool $S$ coincide---is  
\begin{equation}\label{collision_prob}
\Pr_{S}\left(\vct{h}(\tau) = \vct{h}(\tau^*)\right) = \sum_{\mateulerflow} \left( \prod_{r,s=1}^d {{\overlap_{rs}}\choose{\eulerflow_{rs}}} \probquery^{\eulerflow_{rs}} (1-\probquery)^{\overlap_{rs} - \eulerflow_{rs}}\right) \indi \left\{\sum_{s=1}^d \eulerflow_{rs} = \sum_{s=1}^d \eulerflow_{sr} ~,~ \forall r \in [d] \right\},
\end{equation}
where the outer sum is over all arrays of integer numbers $\mateulerflow = (\eulerflow_{rs})_{1 \le r,s\le d}$ such that $0 \le \eulerflow_{rs} \le \overlap_{rs}$ for all $r,s$. We see from the above expression that the collision probability of $\tau$ and $\tau^*$ only depends on the overlap matrix $\matoverlap(\tau, \tau^*)$. We henceforth denote the probability in equation~\eqref{collision_prob} by $q(\matoverlap)$, where we dropped the dependency on $\tau$ and $\tau^*$. Remark that $\tau = \tau^*$ if and only if their overlap matrix $\matoverlap$ is diagonal. Thus, we can rewrite the expected number of solutions as
\begin{equation}
\E[\Zpart-1] = {{\numvar}\choose{\numvar\vct{\pi}}}^{-1} \cdot \sum_{\matoverlap}  {{\numvar}\choose{\matoverlap}} q(\matoverlap)^\numobs ~ \indi\left\{\sum_{r=1}^\numtype \overlap_{rs} = n\pi_s, ~s\in \{1,\cdots, \numtype\} \right\},
\label{expected}
\end{equation}
where the sum is over all non-diagonal arrays $\matoverlap = (\overlap_{rs})_{1 \le r,s\le d}$ with non-negative integer entries that sum to $n$, and ${{n}\choose{\matoverlap}} = \frac{n!}{\prod_{r,s} \overlap_{rs}!}$. 

\paragraph{The rest of the proof} From here, the proof of Theorem~\ref{free_energy_value} roughly breaks into three parts:

$(i)$ One needs to have tight asymptotic control on the collision probability $q(\matoverlap)$ when any subset of the entries of $\matoverlap$ becomes large. This will be achieved via the Laplace method (see, e.g., \cite{debruijn1970asymptotic}). The outcome of this analysis is an asymptotic estimate that exhibits two different speeds of decay, polynomial or exponential, depending on the ``balancedness" of $\matoverlap$ as its entries become large. This notion of balancedness, namely that $\matoverlap$ must have equal row- and column-sums\footnote{These are exactly the constraints on $\mateulerflow$ showing up in~\eqref{collision_prob}.}, is specific to the histogram setting and departs from the usual ``double stochasticity" that arises in other more classical problems such as Graph Coloring, and Community Detection under the stochastic block model~\cite{achlioptas2004two_journal,achlioptas2004chromatic,coja2016chromatic,bapst2016condensation,banks2016information}.    
As we will explain in the next section, configurations $(\tau,\tau^*)$ with an unbalanced overlap matrix have an exponentially decaying collision probability, i.e., they exhibit weak correlation, and disappear very early on as $\numvar \to \infty$ under the scaling $\numobs =  \gamma \frac{\numvar}{\log \numvar}$. On the other hand, those configurations with balanced overlap exhibit a slow decay of correlation: their collision probability decays only polynomially, and these are the last surviving configurations in expression~\eqref{expected} as $\numvar \to \infty$. 

$(ii)$ Understanding the above-mentioned polynomial decay of $q(\matoverlap)$ requires the evaluation of a multivariate Gaussian integral (which is a product of the above analysis) over the space of constraints of the array $\mateulerflow$ in~\eqref{collision_prob}; the latter being the space of \emph{Eulerian flows} on the graph on $d$ vertices whose edges are weighted by the (large) entries of $\matoverlap$. We show that this integral, properly normalized, evaluates to \emph{the inverse square root of the spanning tree (or forest) polynomial} of this graph. This identity seems to be new, to the best of our knowledge, and may be of independent interest. We therefore provide two very different proofs of it, each highlighting different combinatorial aspects. 

$(iii)$ Lastly, armed with these estimates, we show the existence of, and compute the exact value of, the annealed free energy of the model in the thermodynamic limit, thereby completing the proof of Theorem~\ref{free_energy_value}. This last part requires the analysis of a certain optimization problem involving an entropy term and an ``energy" term accounting for the correlations discussed above. Here we can exactly characterize the maximizing configurations for large $\numvar$, and this allows the computation of the value of $\fenergy(\gamma)$. We note once more that this situation contrasts with the more traditional case of Graph Coloring, where we lack a rigorous understanding of the maximizing configurations of the second moment, except when certain additional constraints are imposed on their overlap matrix.

\section{Bounding the collision probabilities}
Here we provide tight asymptotic bounds on the collision probabilities $q(\matoverlap)$ defined in~\eqref{collision_prob}. 
Consider the following subspace of $\R^{d \times d}$, which will play a key role in the analysis: 
\begin{align}\label{def_Fspace}
\Fspace := \left\{\x \in \R^{d \times d} ~:~ \sum_{s=1}^d x_{rs} = \sum_{s=1}^d x_{sr} ~,~ \forall r \in \{1,\cdots,d\} \right\}.
\end{align}
This is a linear subspace of dimension $(d-1)^2+d$ in $\R^{d \times d}$. 
For $p,q \in (0,1)$, let $\kull{p}{q} = p\log(p/q) + (1-p)\log((1-p)/(1-q))$ be the Kullback-Leibler divergence.
Let $G = (V,E)$ be an undirected graph on $d$ vertices where we allow up to two parallel edges between each pair of vertices, i.e., $V = \{1,\cdots,d\}$, and $E \subseteq \{(r,s) ~:~ r,s \in V,~r \neq s\}$. For $\mateulerflow,\matoverlap \in \Rp^{d \times d}$, $\x \in [0,1]^{d\times d}$ let 
\begin{align}\label{varphi_function}
\varphi_{\matoverlap}(\x) := \sum_{(r,s)\in E} \overlap_{rs}\kull{x_{rs}}{\alpha}.
\end{align}
and recalling that $\odot$ represents the Hadamard product, we let
\begin{align}\label{varphi_minimization}
\vartheta(\mateulerflow,\matoverlap) :=  \underset{\substack{\x \in [0,1]^{d \times d}\\ \bm{M}_G(\x\odot \matoverlap , \mateulerflow) \in \Fspace}}{\min} ~\sum_{(r,s) \in E} \overlap_{rs} \kull{x_{rs}}{\alpha},
\end{align}
where for two $d \times d$ matrices $\mtx{a}, \mtx{b}$, $\bm{M}_G(\mtx{a} , \mtx{b})$ is the $d \times d$ matrix with entries $a_{rs}$ if $(r,s)\in E$ and $b_{rs}$ otherwise.
By strong duality (see, e.g., \cite{boyd02,rockafellar70}), the function~\eqref{varphi_minimization} can be written in the more transparent form  
\begin{align*} 
\vartheta(\mateulerflow,\matoverlap) &= \sup_{\vct{\lambda} \in \R^d} \left\{ \sum_{(r,s)\notin E} \eulerflow_{rs}(\lambda_{r} - \lambda_{s}) + \sum_{(r,s)\in E} \overlap_{rs}\log\left(\frac{e^{\lambda_{r} - \lambda_{s}}}{\alpha + (1-\alpha)e^{\lambda_{r} - \lambda_{s}}}\right) \right\},\\
&= \phi_{\matoverlap}^*(\mtx{\eulerflow}\one - \mtx{\eulerflow}^\intercal\one), 
\end{align*}
where $\phi_{\matoverlap}^*$ is the Legendre-Fenchel transform of the (convex) function 
\[\phi_{\matoverlap}(\vct{\lambda}) :=  - \sum_{(r,s)\in E} \overlap_{rs}\log\left(\frac{e^{\lambda_{r} - \lambda_{s}}}{\alpha + (1-\alpha)e^{\lambda_{r} - \lambda_{s}}}\right).\]
We may note that since $\phi^*_{\matoverlap}$ is convex on $\R^d$, $\vartheta$ is a continuous function of its first argument.
Before we state our bounds on the collision probability, we recall the following concept from algebraic graph theory. Define \emph{the spanning tree polynomial} of $G$ as
\[T_G(\z) := \frac{1}{\nst(G)} ~\sum_{T} \prod_{(r,s) \in T} z_{rs},\]
for $\z \in \Rp^{d \times d}$, where the sum is over all spanning trees of $G$, and $\nst(G)$ is the number of spanning trees of $G$.
In cases where $G$ is not connected, we define the following polynomial 
 \[P_G := \overset{\ncc(G)}{\underset{l=1}{\prod}} T_{G_l},\] 
where $G_l$ is the $l$th connected component of $G$, and we denote by $\ncc(G)$ the number of connected components of $G$. This polynomial may be interpreted as the generating polynomial of \emph{spanning forests} of $G$ having exactly $\ncc(G)$ trees. The polynomials $T_{G}$ and $P_{G}$ are multi-affine, homogenous of degree $d-1$ for $T_G$ (when $G$ is connected) and $d-\ncc(G)$ for $P_{G}$, and do not depend on the diagonal entries $\{z_{rr} : 1 \le r \le d\}$. Furthermore, letting $z_{rs} = 1$ for all $r\neq s$, we have $P_{G}(\z) = T_G(\z) = 1$.
We now provide tight asymptotic bounds on the collision probability $q(\matoverlap)$ when a subset $E$ of the entries of $\matoverlap$ become large.
\begin{theorem} \label{collision_asymptotic}
 Let $G = (V,E)$ with $V=\{1,\cdots,d\}$, $E = \{(r,s) \in V^2~:~ r\neq s\}$, and $\epsilon \in (0,1)$. There exist two constants $0 <c_u < c_l$ depending on $\epsilon$, $d$ and $\alpha$ such that for all $n$ sufficiently large, and all $\matoverlap \in \{0,\cdots,n\}^{d \times d}$ with $\overlap_{rs} \ge \epsilon n$ if and only if $(r,s) \in E$, we have
\begin{align*}
 c_l \frac{e^{-\vartheta_l(\matoverlap)}}{P_{G}(\matoverlap)^{1/2}} \le 
 q(\matoverlap) \le c_u \frac{e^{-\vartheta_u(\matoverlap)}}{P_{G}(\matoverlap)^{1/2}}.
\end{align*}
with 
\[\vartheta_u(\matoverlap)  = \inf_{\mateulerflow} \{\vartheta(\mateulerflow,\matoverlap) : 0 \le \eulerflow_{rs}\le \overlap_{rs}~ \forall (r,s) \notin E \},\]
and 
\[\vartheta_l(\matoverlap)  = \sup_{\mateulerflow} \{\vartheta(\mateulerflow,\matoverlap) : 0 \le \eulerflow_{rs}\le \overlap_{rs}~ \forall (r,s) \notin E \}.\]\end{theorem} 
  
Let us now expand on the above result and derive some special cases and corollaries. First, we see that the collision probabilities can decay at two different speeds---polynomial or exponential---in the entries of the overlap matrix $\matoverlap$, depending on whether $\vartheta_u(\matoverlap)$ (and/or $\vartheta_l(\matoverlap)$) is zero or strictly negative. Second, the apparent gap in the exponential decay of $q(\matoverlap)$ in the above characterization is artificial; one can make $\vartheta_u$ and $\vartheta_l$ equal by taking $\overlap_{rs} = 0$ for all $(r,s) \notin E$.  Alternatively, they could be made arbitrarily close to each other under an appropriate limit: Assume for simplicity that $\overlap_{rs} = n w_{rs}>0$ for all $(r,s) \in E$ for some $\w \in [0,1]^{d \times d}$. We have
\[\vartheta(\mateulerflow,\matoverlap) = n \vartheta(\mateulerflow/n,\w).\]
For $(r,s) \notin E$, we have $\overlap_{rs} < \epsilon n$, therefore
\[\vartheta_u(\matoverlap)/n \le \inf_{\x} \{\vartheta(\x,\w) : 0 \le x_{rs} \le \epsilon~ \forall (r,s) \notin E \} \underset{\epsilon \to 0}{\longrightarrow} \vartheta(\vct{0},\w).\] 
The last step is justified by the continuity of $\vartheta(\cdot,\w)$. The same argument holds for $v_l(\matoverlap)$. Denoting the limiting function under this operation as $\vartheta(\w)$, we obtain:
\[\vartheta(\w) = \sup_{\vct{\lambda} \in \R^d} \sum_{(r,s)\in E} w_{rs}\log\left(\frac{e^{\lambda_{r} - \lambda_{s}}}{\alpha + (1-\alpha)e^{\lambda_{r} - \lambda_{s}}}\right) = \underset{\substack{\x \in [0,1]^{d \times d}\\ \w \odot \x\in \Fspace}}{\min}~\varphi_{\w}(\x) .\]

The function $\vartheta$ can be seen as the exponential rate of decay of $q(\matoverlap)$. The reason $\vartheta_u$ and $\vartheta_l$ cannot (in general) be replaced by $\vartheta$ in Theorem~\ref{collision_asymptotic} is that all control on the constants $c_u$ and $c_l$ is lost when $\epsilon \to 0$.
 Next, we identify the cases where this exponential decay is non-vacuous.
\begin{lemma}\label{inFornot}
Let $\alpha \in (0,1)$, and $\matoverlap \in \Rp^{d \times d}$. We have 
\begin{itemize}
\item[(i)] $\vartheta(\matoverlap) = 0$ if and only if $\matoverlap \in \Fspace$, 
\item[(ii)] $\vartheta_u(\matoverlap) = 0$ if and only if $\bm{M}_G(\alpha\matoverlap,\mateulerflow) \in \Fspace$ for \emph{some} $\mateulerflow \in \Rp^{d\times d}$ such that $0 \le \eulerflow_{rs} \le \overlap_{rs}$ for all $(r,s) \notin E$.

\end{itemize}
\end{lemma}

Now we specialize Theorem~\ref{collision_asymptotic} to the case where the entries of the overlap matrix are either zero or grow proportionally to $n$. From Theorem~\ref{collision_asymptotic} and Lemma~\ref{inFornot}, we deduce a key corollary on the convergence of the properly rescaled logarithm of the collision probabilities.
\begin{corollary}\label{free_energy_collision}
 Given a graph $G = (V,E)$, let $\w \in [0,1]^{d \times d}$ be such that $w_{rs} >0$ if and only if $(r,s) \in E$.
 \noindent If $\w \in \Fspace$ then 
 \[\lim_{n \to \infty} \frac{\log q(n\w)}{\log n} = - \frac{d-\ncc(G)}{2}.\]
Otherwise if $\w \notin \Fspace$, then 
 \[\lim_{n \to \infty} \frac{\log q(n\w)}{n} = - \vartheta(\w).\]
\end{corollary}
We see that the assignments $\tau$ such that $\matoverlap(\tau,\tau^*)\in \Fspace $ exhibit a much stronger correlation to $\tau^*$ than those for which this overlap matrix does not belong to $\Fspace$, and will hence survive much longer as $n \to \infty$.   
    
\vspace{.5cm}
\begin{proofof}{Lemma~\ref{inFornot}} 
Let $\matoverlap,\mateulerflow \in \R_+^{d \times d}$ with $\matoverlap \neq \mtx{0}$. Let $\alpha \in (0,1)$, and let $G=(V,E)$ denote a graph on $d$ vertices. 
The function $\varphi_{\matoverlap}$ defined in~\eqref{varphi_function} is strictly convex on the support of $\matoverlap$, i.e., on the subspace induced by the non-zero coordinates of $\matoverlap$, so it admits a unique minimizer on the closed convex set $ \{\x \in [0,1]^{d \times d} ~:~ \bm{M}_{G}(\x^*\odot \matoverlap,\mateulerflow) \in \Fspace\}$ intersected with that subspace. Let $\x^*$ be this minimizer. 
By differentiating the associated Lagrangian, the entries of $\x^*$ admit the expressions
\[x^*_{rs} = \frac{\alpha}{\alpha +(1-\alpha)e^{\lambda_{r} - \lambda_{s}}},\]  
for all $ (r,s) \in E$ (recall that $\overlap_{rs} >0$ for all such $(r,s)$), and where the vector $\vct{\lambda} \in \R^d$ is the unique solution up to global shifts of the system of equations 
\begin{equation}\label{equation_lambda}
\sum_{s: (r,s)\in E} \frac{\alpha\overlap_{rs}}{\alpha+(1-\alpha)e^{\lambda_r-\lambda_s}} + \sum_{s:(r,s)\notin E} \eulerflow_{rs}= \sum_{s: (r,s)\in E} \frac{\alpha\overlap_{sr}}{\alpha+(1-\alpha)e^{\lambda_s-\lambda_r}} + \sum_{s:(r,s)\notin E} \eulerflow_{sr} \quad \forall r \in \{1,\cdots,d\}.
\end{equation}
The claims of the lemma follow directly from the system of equations~\eqref{equation_lambda} and the fact that the non-negative function $\varphi_{\matoverlap}$ vanishes if and only if $x^*_{rs}= \alpha$ for all $(r,s)\in E$: to show $(i)$, we take $\mateulerflow = \mtx{0}$. It is clear from the equations that $\matoverlap \in \Fspace$ if and only if $\vct{\lambda} = c\one$, $c \in \R$, is a solution to the above equations; and this is equivalent to $x^*_{rs} = \alpha$ whenever $\overlap_{rs} >0$. This is in turn equivalent to $\vartheta(\matoverlap) = \varphi_{\matoverlap}(\x^*) = 0$. The same strategy is employed to show $(ii)$, in conjunction with the continuity of the function $\mateulerflow \mapsto \vartheta(\mateulerflow,\matoverlap)$ over a compact domain (the infimum defining $\vartheta_u$ is attained). 
\end{proofof}

\vspace{.5cm}
\begin{proofof}{Corollary~\ref{free_energy_collision}}
 Fix $G=(V,E)$, let $\w \in (0,1)^{d\times d}$ with $w_{rs}>0$ if and only if $(r,s)\in E$, and let $n$ be an integer. For simplicity, assume that for $n\w$ is an array of integer entries. The non-integer part introduces easily manageable error terms. 
Applying Theorem~\ref{collision_asymptotic} with $\epsilon = \min_{(r,s)\in E}w_{rs}$, we have for $n$ large
 \[c_l P_G(n\w)^{-1/2} \exp-\vartheta_l(n\w) \le q(n\w) \le c_u P_G(n\w)^{-1/2} \exp-\vartheta_u(n\w).\] 
Moreover, since $w_{rs}=0$ for $(r,s)\notin E$, we have   
 \[\vartheta_u(n\w) = \vartheta_l(n\w) = n \vartheta(\w).\]
On the other hand, by homogeneity of the polynomial $P_G$, $P_G(n\w) = n^{d-\ncc(G)}P_{G}(\w)$. 
Applying Lemma~\ref{inFornot} yields the desired result: If $\w \in \Fspace$ then
 \[\lim_{n \to \infty} \frac{\log q( n \w)}{\log n} = - \frac{d-\ncc(G)}{2}.\]
Otherwise,
 \[\lim_{n \to \infty} \frac{\log q( n \w)}{n} = - \vartheta(\w).\]
\end{proofof}

\subsection{A Gaussian integral}
One important step in proving Theorem~\ref{collision_asymptotic} (specifically for obtaining the polynomial decay part of $q(\matoverlap)$) is the following identity relating the Gaussian integral on a linear space $\Fspace(G)$ defined based on a graph $G$ to the spanning tree/forest polynomial of $G$.  We denote by $K_d$ the complete graph on $d$ vertices where every pair of distinct vertices is connected by \emph{two parallel edges}.  
\begin{proposition}\label{gaussian-integral}
Let $G = (V,E)$ be a graph on $d$ vertices, where self-loops and up to two parallel edges are allowed: $V = \{1,\cdots,d\}$, $E \subseteq V \times V$. Further, let 
\[\Fspace(G) = \Big\{\x \in \Fspace ~:~ x_{rs} = 0 ~\mbox{for}~(r,s) \notin E\Big\}.\] 
For any array of positive real numbers $(w_{rs})_{(r,s) \in E}$, we have  
\[\int_{\Fspace(G)} e^{-\sum_{rs} x_{rs}^2/2w_{rs}} ~\mathrm{d}\x = \left((2\pi)^{\dim(\Fspace(G))}~\frac{\prod_{r,s} w_{rs}}{P_G(\w)}\right)^{1/2}.\]
\end{proposition} 
In the case where $G$ is the complete graph $K_d$, $\Fspace(G) = \Fspace$, $\dim(\Fspace) = (d-1)^2+d$, and $P_G = T_G = (2^{d-1}d^{d-2})^{-1} \sum_{T} \prod_{(r,s)\in T}w_{rs}$ where the sum is over all spanning trees of $K_d$.  The pre-factor in the last expression comes from Cayley's formula for the number of spanning trees of the complete graph. We will show that it suffices to prove Proposition~\ref{gaussian-integral} in the case where $G= K_d$ in order to establish it for any graph $G$. We were not able to locate this identity in the literature. To illuminate the combinatorial mechanisms behind it, we provide what appear to be two very different proofs of it. A first ``direct" and purely combinatorial proof views $\Fspace(G)$ as the space of \emph{Eulerian flows} of the  graph $G$. A second, slightly indirect proof which is mainly analytic, and relates the above Gaussian integral to the characteristic polynomial of the Laplacian matrix of $G$ then invokes the Principal Minors Matrix Tree theorem (see, e.g., \cite{chaiken1982combinatorial}).  

\section{Computing the annealed free energy}

In this section we establish the existence of $\fenergy(\gamma)$, and compute its value for all $\gamma >0$.   
For $1\le k \le d$ let $\mathcal{D}_k$ denote the set of binary matrices $\bm{X} \in \{0,1\}^{k\times d}$ such that each column of $\bm{X}$ contains \emph{exactly} one non-zero entry and each row contains \emph{at least} one non-zero entry. The elements of $\mathcal{D}_k$ represent partitions of the set $\{1,\cdots,d\}$ into $k$ non-empty subsets.
\begin{proposition}\label{free_energy}
Let $m = \gamma \frac{n}{\log n}$ with $\gamma >0$ fixed for all $n \ge 2$. We have 
\begin{align*}
\fenergy (\gamma) =\max_{ 1 \le k \le d-1} \left\{ H(\vct{\pi}) - \min_{\bm{X} \in \mathcal{D}_k} H(\bm{X}\vct{\pi}) - \frac{\gamma}{2}(d-k)\right\}.
\end{align*}
\end{proposition}
\noindent Moreover, the inner minimization problem in the above expression can be solved explicitly:
\begin{lemma}\label{nncvx_maximization}
Let $\vct{\pi}_{[\cdot]}$ be a permutation of the vector $\vct{\pi}$ such that $\pi_{[1]} \ge \pi_{[2]} \ge \cdots \ge \pi_{[d]}$. And for $1\le k \le d-1$, let $\vct{\pi}^{(k)} \in \Delta^{k-1}$ defined as 
${\pi}^{(k)}_{1} = \sum_{r=1}^{d-k+1} \pi_{[r]}$ and ${\pi}^{(k)}_{l} = \pi_{[d-k+l]}$ for all $2\le l \le k$ (if $k\ge 2$).
Then
\[\min_{\bm{X} \in \mathcal{D}_k}H(\bm{X}\vct{\pi}) = H(\vct{\pi}^{(k)}).\]
\end{lemma}
\noindent Theorem~\ref{free_energy_value} follows from Proposition~\ref{free_energy} and Lemma~\ref{nncvx_maximization}. We begin with the proof of the latter and devote the next subsection to the lengthier proof of the former.

\vspace{0.1in}
\noindent \begin{proofof}{Lemma~\ref{nncvx_maximization}}
We start with an arbitrary partition of $\vct{\pi}$ into $k$ groups, and define a sequence of operations on the set of $k$-partitions of $\vct{\pi}$ that strictly decreases $H(\mtx{X}\vct{\pi})$ at each step, and, irrespective of the starting point, always converges to $\vct{\pi}^{(k)}$. Starting with an arbitrary $k$-partition, write down the groups from left to right in decreasing order of total weight of each group. Initially, every group is marked \emph{incomplete}. Then we perform the following operations:
\begin{enumerate}
\item Start with the rightmost incomplete group. 
\item If it has more than one element, transfer the largest element to the leftmost group. This strictly decreases the entropy, since the heaviest group gets heavier and the lightest group gets lighter. Repeat this step until the rightmost group has exactly one element, and then move to the next step.
\item Consider this (now singleton) group. If there is no element to its left that is lighter than it, mark the group as complete. Else, swap this element with the lightest element to its left, and then mark it complete. Then go back to step 1.
\end{enumerate}
\end{proofof}

\subsection{Proof of Proposition~\ref{free_energy}}
Let $\numobs = \gamma \frac{n}{\log n}$. Recall from equation~\eqref{expected} that 
\[\E[\Zpart-1] = {{\numvar}\choose{\numvar\vct{\pi}}}^{-1} \cdot \sum_{\matoverlap}  {{\numvar}\choose{\matoverlap}} q(\matoverlap)^{\numobs} ~ \indi\left\{ \matoverlap^\intercal \one = n \vct{\pi} \right\},\]
where the sum is over all arrays $\matoverlap \in \Zp^{d\times d}$ such that $\one^\intercal \matoverlap \one = n$, $1\le \sum_{r\neq s} \overlap_{rs}$.
Since the sum defining $\E[\Zpart-1]$ is larger than its maximum term and smaller than the maximum term times $(n+1)^{d^2}$, we only need to understand the convergence of the sequence  
\begin{align*} 
\fenergy_n &:= \frac{1}{n} \log \left( \max_{\substack{\matoverlap \in \{0,\cdots,n\}^{d \times d}  \\ \text{ non-diagonal}}}{{n}\choose{\matoverlap}} q(\matoverlap)^{\numobs} ~\indi\left\{ \matoverlap^\intercal \one = n \vct{\pi} \right\} \right)\\
&= \max \left\{ \frac{1}{n} \log {{n}\choose{ \matoverlap }} + \gamma \frac{\log q( \matoverlap )}{\log n} ~:~ \matoverlap \in \{0,\cdots,n\}^{d \times d}, \sum_{r\neq s} \overlap_{rs} \ge 1, \matoverlap^\intercal \one = n \vct{\pi} \right\}.
\end{align*}
If this sequence converges, we would have
\begin{equation}\label{free_limit_1}
\fenergy(\gamma) = -H(\vct{\pi}) + \lim_{n \to \infty} \fenergy_n,
\end{equation}
since $\frac{1}{n} \log {{n}\choose{ n\vct{\pi} }} \to H(\vct{\pi})$ by Stirling's formula. Next, we show that the above limit indeed exists. 
Let 
\begin{align}\label{psi_n_function}
\psi_n(\w) := \frac{1}{n} \log {{n}\choose{ n \w }} + \gamma \frac{\log q( n \w )}{\log n}.
\end{align}
By Corollary~\ref{free_energy_collision}, the function 
\begin{align}\label{psi_function}
\psi(\w) := 
\begin{cases}
H(\w) - \frac{\gamma}{2}(d-\ncc(\w)) & \mbox{if} ~ \w \in \Fspace, \\
-\infty & \mbox{otherwise},
\end{cases}
\end{align}
is the point-wise limit of the sequence of functions $\{\psi_n\}_{n \ge 2}$ on $\Delta^{d \times d-1}$.
Next, we use the following lemma which states that any non-diagonal sequence of maximizers $\{\matoverlap^{(n)}\}_n$ of $\psi_n$ is such that $\sum_{r\neq s}\overlap^{(n)}_{rs}$ grows proportionally to $n$.
\begin{lemma}\label{maximizing_psi}
For all $n\ge 2$, let 
\[\matoverlap^{(n)} \in \arg\max \left\{\psi_{n}(\matoverlap/n) ~:~ \matoverlap \in \{0,\cdots,n\}^{d \times d}, ~ 1 \le \sum_{r\neq s} \overlap_{rs} \le n, ~\matoverlap^\intercal \one= n \vct{\pi} \right\}.\]
It holds that
\[\liminf_{n\to \infty} ~ \frac{\sum_{r\neq s}\overlap^{(n)}_{rs}}{n} >0. \] 
\end{lemma}
By Lemma~\ref{maximizing_psi}, which we prove at the end of the current argument, we can safely restrict the set of candidate maximizers to those $\matoverlap$ such that $\sum_{r\neq s} \overlap_{rs} \ge  c_0 n$ for some fixed but small $c_0>0$. 
From here, and by a change of variables $\matoverlap = n \w$, mere point-wise convergence suffices to interchange $\liminf$ and $\sup$:   
\begin{align}
\underset{n \to \infty}{\lim\inf} ~\fenergy_n
&\ge \underset{n \to \infty}{\lim\inf} ~\sup \left\{ \psi_n(\w) ~:~ \w \in \{i/n : 0 \le i \le n\}^{d \times d},~ c_0 \le \sum_{r\neq s} w_{rs} \le 1, ~\w^\intercal \one = \vct{\pi} \right\}\nonumber\\
&\ge \sup \left\{ \psi(\w) ~:~ \w \in [0,1]^{d \times d} \cap \Fspace, ~c_0 \le \sum_{r\neq s} w_{rs} \le 1, ~\w^\intercal \one = \vct{\pi} \right\}.
\label{liminf}
\end{align}

\vspace{.3cm}

Now we present a matching upper bound for $\limsup \fenergy_n$.  
For $\epsilon >0$, let $G_n = (\{1,\cdots,d\},E_n)$ be defined such that $(r,s) \in E_n$ if and only if $w^{(n)}_{rs}\ge \epsilon$. Let $(G_l)_{l=1}^k$ denote the connected components of the graph $G_n$, $k = \ncc(G_n)$. Also, for $\w$ an array for positive entries, let $\ncc^\epsilon(\w)$ denote the number of connected components of the graph $G(\w,\epsilon)=(V,E(\w,\epsilon))$, $V = \{1,\cdots,d\}$, $E(\w,\epsilon)=\{(r,s)~:~ r\neq s, ~ w_{rs} > \epsilon\}$, and let
\[\vartheta^\epsilon(\w) := \inf_{\x} \{\vartheta(\x,\w) : 0 \le x_{rs} \le \epsilon~ \forall (r,s) \notin E(\w,\epsilon)\}.\] 
We will also write $\ncc(\w)$ for $\ncc^0(\w)$.
Let $\w^{(n)} = \matoverlap^{(n)}/n$ for all $n \ge 2$, where $\matoverlap^{(n)}$ is defined in Lemma~\ref{maximizing_psi}. By Theorem~\ref{collision_asymptotic}, we have for $n$ sufficiently large
\[q(n\w^{(n)}) \le c_u(\epsilon,d,\alpha) P_{G_n}(n\w^{(n)})^{-1/2} \exp - \vartheta^{\epsilon} (n\w^{(n)}).\]
Since $w^{(n)}_{rs} \ge \epsilon$ of all the edges $(r,s)$ of $G_n^\epsilon$, $\prod_l T_{G_l}(\w^{(n)})$ is bounded below by $\epsilon^{d}$ \emph{independently of $n$}. 
Therefore, for $n$ sufficiently large, 
\begin{align*}
\psi_n(\w^{(n)}) &=  \frac{1}{n} \log {{n}\choose{n\w^{(n)}}} + \gamma \frac{\log q(n\w^{(n)})}{\log n} \\
&\le \frac{1}{n} \log {{n}\choose{n\w^{(n)}}} - \frac{\gamma}{2}(d - \ncc^{\epsilon}(\w^{(n)})) - \frac{\gamma n}{\log n} \vartheta^{\epsilon}(\w^{(n)})+ \bigo \left(\frac{\log c_u(\epsilon,d,\alpha)+d\log(1/\epsilon)}{\log n}\right)\\
&\le \sup \left\{ H(\w) - \frac{\gamma}{2}(d - \ncc^{\epsilon}(\w)) - \frac{\gamma n}{\log n} \vartheta^{\epsilon}(\w)~:~ \w \in [0,1]^{d\times d}, c_0 \le \sum_{r\neq s} w_{rs} \le 1, ~\w^\intercal \one = \vct{\pi} \right\} \\
& \quad + \bigo \left(\frac{\log c_u(\epsilon,d,\alpha)+d\log(1/\epsilon)}{\log n}\right),
\end{align*}
where the last inequality is obtained by Stirling's formula and taking a supremum over all $\w$. By Lemma~\ref{inFornot}, $\vartheta^{\epsilon}(\w) = 0$ if and only if $\bm{M}_{G}(\alpha\w,\x) \in \Fspace$ for some $\x \in [0,1]^{d\times d}$ such that $0 \le x_{rs} \le \epsilon$ for all $(r,s) \notin E$, $G=(V,E)$ being the graph whose edges are $(r,s) : w_{rs} \ge \epsilon$. This constrains the supremum to be achieved in the space of such $\w$ for $n$ sufficiently large. Moreover, this condition implies in particular that 
\[\|\w \one - \w^\intercal \one\|_{\ell_{\infty}} \le 2d \alpha^{-1}\epsilon,\]
where $\|\cdot \|_{\ell_{\infty}}$ is the $\ell_{\infty}$ norm of a vector in $\R^{d}$.
Consequently, this yields the following upper bound as $n \to \infty$,
\begin{align}\label{limsup}
\underset{n \to \infty}{\limsup} ~\fenergy_n \le 
\sup \left\{ 
H(\w) - \frac{\gamma}{2}(d - \ncc^{\epsilon}(\w))
~:~ 
\begin{array}{c}
\w \in [0,1]^{d\times d}, \|\w \one - \w^\intercal \one\|_{\ell_{\infty}} \le 2d \alpha^{-1}\epsilon, \\
c_0 \le \sum_{r\neq s} w_{rs} \le 1, ~\w^\intercal \one = \vct{\pi}
\end{array}
\right\},
\end{align}
for all $\epsilon>0$. Next, we argue that as $\epsilon \to 0$, the right-hand side of the above inequality converges to 
\[ \sup \left\{ H(\w) - \frac{\gamma}{2}(d - \ncc(\w))~:~ \w \in [0,1]^{d\times d}\cap \Fspace, c_0 \le \sum_{r\neq s} w_{rs} \le 1,~\w^\intercal \one = \vct{\pi} \right\},\]
thereby establishing the existence of the limit $\lim \fenergy_n$ along with its precise value.
Since the function $\epsilon \to \ncc^\epsilon(\w)$ is non-decreasing for any fixed $\w$, the limit of the right-hand side of~\eqref{limsup} as $\epsilon \to 0$ exists by monotone convergence. The limit can be decomposed as 
\begin{align*} 
&\lim_{\epsilon \to 0} ~\sup \left\{ H(\w) - \frac{\gamma}{2}(d - \ncc^{\epsilon}(\w))
~:~ 
\begin{array}{c}
\w \in [0,1]^{d\times d}, \|\w \one - \w^\intercal \one\|_{\ell_{\infty}} \le 2d \alpha^{-1}\epsilon, \\
c_0 \le \sum_{r\neq s} w_{rs} \le 1, ~\w^\intercal \one = \vct{\pi}
\end{array}
 \right\}\\
&= \max_{1\le k \le d} \max_{\{V_l\}_{l=1}^k} \lim_{\epsilon \to 0} ~ \sup 
\left\{
H(\w) - \frac{\gamma}{2}(d - k)~:
\begin{array}{c}
\w \in [0,1]^{d \times d}, ~\|\w \one - \w^\intercal \one\|_{\ell_{\infty}} \le 2d \alpha^{-1}\epsilon, \\
w_{rs} \le \epsilon ~\forall (r,s) \in V_l \times V_{l'}, ~ l\neq l',\\
G_l(\w) \text{ is connected}~\forall l, ~ c_0 \le \sum_{r\neq s} w_{rs} \le 1, ~ \w^\intercal \one = \vct{\pi}
\end{array}
\right\},
\end{align*}
where $\{V_l\}_{l=1}^k$ ranges over a partitions of the set $\{1,\cdots,d\}$ with $k$ non-empty subsets, and $G_l(\w) = (V_l, \{(r,s)\in V_l \times V_l~:~ w_{rs}>\epsilon\})$ for all $1\le l \le k$. Letting $\epsilon < c_0$, the range of the outer-most maximum becomes $1\le k \le d-1$.
By concavity of the entropy, the constraint that the graphs $G_l(\w)$ must be connected can be safely removed from the maximization problem without changing its maximum value since it will be automatically satisfied. Thus, the inner-most optimization problem is that of a continuous function on a closed and bounded domain that shrinks with $\epsilon$. Its value is therefore a continuous function of $\epsilon$. Hence, by sending $\epsilon$ to $0$, in conjunction with the lower bound~\eqref{liminf}, we conclude that
\begin{equation}\label{limit_free}
\lim_{n\to \infty} \fenergy_n = \sup \left\{ \psi(\w) ~:~ \w \in [0,1]^{d \times d}, ~c_0 \le \sum_{r\neq s} w_{rs} \le 1, ~\w \one = \w^\intercal \one = \vct{\pi} \right\}.
\end{equation} 

As a final step, we make the above expression a bit more explicit. As argued previously, the supremum in~\eqref{limit_free} can be decomposed such that one first takes the maximum of $\psi(\w)$ over all $\w$ such that $w_{rs} = 0$ for all $(r,s) \in V_l \times V_{l'}, ~ l \neq l'$ where $\{V_l\}_{1\le l \le k}$ is a fixed partition of $\{1,\cdots,d\}$ into non-empty subsets, then maximize over all such partitions, then over all $1\le k\le d-1$. The first optimization problem has a value  
\begin{align*}
\sup \left\{ H(\w) - \frac{\gamma}{2} (d-k) ~:~ \w \in [0,1]^{d \times d},~ w_{rs} = 0, ~ (r,s) \in V_l \times V_{l'}, ~ l \neq l', ~ \w \one = \w^\intercal \one = \vct{\pi} \right\},
\end{align*}
where the constraint $c_0 \le \sum_{r\neq s} w_{rs} \le 1$ is not active for $c_0$ small enough, hence can be removed. Let $\w$ be in the above constraint set. Then $H(\w) = - \sum_{l=1}^k \sum_{(r,s)\in V_l \times V_l}w_{rs} \log w_{rs}$, and this is maximized at 
\begin{align} \label{w_maximizer}
w^*_{rs} = 
\begin{cases}
(\pi_r \pi_s)/\sum_{r' \in V_l} \pi_{r'} & \mbox{ if } (r,s) \in V_l \times V_l, ~ l \in \{1,\cdots,k\}, \\
0 & \mbox{ otherwise, }
\end{cases}
\end{align}
with maximum value 
\begin{align} \label{maximum_entropy}
H(\w^*) &= 2 H(\vct{\pi}) + \sum_{l=1}^k \left(\sum_{r \in V_l}\pi_{r}\right) \log \left(\sum_{r \in V_l}\pi_{r}\right),\\
&= 2 H(\vct{\pi}) - H(\bm{X}\vct{\pi}),\nonumber
\end{align}
where $\bm{X} \in \{0,1\}^{k\times d}$, $X_{l,r} = 1$ if and only if $r \in V_l$. Note that $\mathcal{D}_k$ is the set of all such matrices (each one corresponding to a  partition $\{V_l\}$ of $\{1,\cdots,d\}$). 
Finally, by maximizing over all possible partitions, and using~\eqref{free_limit_1} we get 
\[\fenergy (\gamma) = \max_{1 \le k \le d-1} \left\{H(\vct{\pi}) - \min_{\bm{X} \in \mathcal{D}_k} H(\bm{X}\vct{\pi}) - \frac{\gamma}{2} (d-k)\right\}.\]  
\noindent This completes the proof of Proposition~\ref{free_energy}, except for the proof Lemma~\ref{maximizing_psi}, which we provide below.

\vspace{.5cm}

\noindent 
\begin{proofof}{Lemma~\ref{maximizing_psi}}
Let 
\[\matoverlap^{(n)} \in \arg\max \left\{\psi_{n}(\matoverlap/n) ~:~ \matoverlap \in \{0,\cdots,n\}^{d \times d}, ~ 1 \le \sum_{r\neq s} \overlap_{rs}, ~\matoverlap^\intercal \one= n \vct{\pi} \right\}.\]
We show that 
\[\liminf_{n\to \infty} ~ n^{-1} \sum_{r\neq s}\overlap^{(n)}_{rs} > 0.\]
Let us first show that  
\[ \frac{(\log n)^3}{n} \sum_{r\neq s}\overlap^{(n)}_{rs} \longrightarrow \infty,\]
and then remove the logarithmic factor.
We proceed by contradiction, by showing that if the above statement is not true, then the expected number of non-planted solutions $\E[\Zpart-1]$ vanishes as $n \to \infty$ for any $\gamma >0$, which contradicts our  lower bound of Theorem~\ref{mainthm}.  We have  
\[\E\left[\Zpart -1\right] \le {{n}\choose{n \vct{\pi}}}^{-1} \cdot  (n+1)^{d^2} \cdot {{n}\choose{\matoverlap^{(n)}}} \cdot q_{\max}^{\gamma n/\log n} ,\]
with $q_{\max} = \max \left\{ q(\matoverlap) ~:~ 1 \le \sum_{r\neq s} \overlap_{rs} \right\} < 1$. Moreover, 
\[{{n}\choose{\matoverlap^{(n)}}} =  {{n}\choose{n \vct{\pi}}} \prod_{r=1}^d \frac{(n\pi_r)!}{\prod_{s \neq r} \overlap_{sr}! (n\pi_r - \sum_{s \neq r}\overlap_{sr})!} \le {{n}\choose{n \vct{\pi}}}  \prod_{r=1}^d (n\pi_r)^{\sum_{s\neq r} \overlap_{sr}}.\]
If $\sum_{r\neq s}\overlap^{(n)}_{rs} \le C n/(\log n)^3$ for some constant $C>0$, then
\[\E\left[\Zpart -1\right] \le (n+1)^{d^2} \cdot n^{Cn/(\log n)^3} \cdot q_{\max}^{\gamma n/\log n} \underset{n \to \infty}{\longrightarrow} 0,\]
for all $\gamma >0$, and this contradicts the fact that below $\gamma_{\text{low}}$ there are exponentially many distinct satisfying assignments.

Now let us assume that $ \frac{(\log n)^3}{n}\sum_{r\neq s}\overlap^{(n)}_{rs} \to \infty$ but $\liminf n^{-1} \sum_{r\neq s}\overlap^{(n)}_{rs} = 0$. We proceed by contradiction once more, and construct a sequence of points that have a higher objective value than $\matoverlap^{(n)}$.   
 Instead of working with convergent subsequences, we may as well assume that $\{\matoverlap^{(n)}\}$ is convergent.    
 Let 
 \[E_n = \left\{(r,s)~:~ r \neq s, ~ \overlap_{rs}^{(n)} > \epsilon \sum_{r \neq s} \overlap_{rs}^{(n)}\right\},
\quad \text{and} \quad
 E_{\infty} = \left\{(r,s) ~:~ r\neq s,~ \liminf_{n\to \infty} \frac{\overlap^{(n)}_{rs}}{\sum_{r\neq s} \overlap^{(n)}_{rs}} > 0\right\},\] 
 for all $n$ and some $\epsilon >0$ sufficiently small.
Let $k_n = \ncc(G_n)$ be the number of connected components of the graph $G_n = (\{1,\cdots, d\},E_n)$, and similarly, let $k_{\infty} = \ncc(G_\infty)$ with $G_\infty = (\{1,\cdots, d\},E_\infty)$. Observe that $E_\infty$ and $E_n$ are both non-empty sets, hence $k_\infty, k_n \le d-1$ for all $n$.

Now we consider an arbitrary partition of the set of vertices $\{1,\cdots,d\}$ into $k_\infty$ subsets $\{V_l\}_{1\le l \le k_\infty}$, and let $G$ be the graph on $d$ vertices with edge set $\cup_{l=1}^{k_{\infty}} V_l \times V_l$; i.e., $G$ is the union of $k_\infty$ \emph{complete} connected components. Finally, let $\mtx{v}^{(n)} := n \w$ for all $n$, with 
\[w_{rs} = 
\begin{cases}
(\pi_r \pi_s)/\sum_{r' \in V_l} \pi_{r'} & \mbox{ if } (r,s) \in V_l \times V_l, ~ l \in \{1,\cdots,k_{\infty}\}, \\
0 & \mbox{ otherwise, }
\end{cases}
\]
Recall that this construction provides one of the candidate maximizers of the annealed free energy (see~\eqref{w_maximizer}). Observe that $\mtx{v}^{(n)}$ satisfies all the constraints satisfied by $\matoverlap^{(n)}$, and additionally, $\mtx{v}^{(n)} \in \Fspace$. Therefore, by Corollary~\ref{free_energy_collision}, we have
\[\psi_n(\mtx{v}^{(n)}/n) = H(\w) -\frac{\gamma}{2} (d-k_{\infty}) + o_n(1).\] 
Recall that the function $\psi_n$ is defined in~\eqref{psi_n_function}. On the other hand, to study the asymptotics of $\psi_n(\matoverlap^{(n)}/n)$, we apply Theorem~\ref{collision_asymptotic} with $n$ replaced by $\sum_{r\neq s}\overlap^{(n)}_{rs}$ (which grows to infinity), and we get
\[\psi_n(\matoverlap^{(n)}/n) \le H(\vct{\pi}) -\frac{\gamma}{2} (d-k_n)\left(1-3\frac{\log\log n}{\log n}\right) - \frac{\vartheta_u(\matoverlap^{(n)})}{\log n} + o_n(1).\]
The term in the right-hand side follows from Stirling's formula and the fact that $\overlap_{rs}^{(n)}/n \to 0$ for all $r \neq s$. The second term follows from the fact that
\[P_{G_n}(\matoverlap^{(n)}) \ge \left(\epsilon\sum_{r\neq s}\overlap^{(n)}_{rs}\right)^{d-k_n} \gg \left(\frac{n}{(\log n)^3}\right)^{d-k_n}. \]
 
 Next, we argue based on these estimates that $\psi_n(\mtx{v}^{(n)}/n) > \psi_n(\matoverlap^{(n)}/n)$ for all $n$ large enough. First, the term involving $\vartheta_u$ in the upper bound on $\psi_n(\matoverlap^{(n)}/n)$ can be dropped since it is always non-negative.
By direct computation (we already showed this in~\eqref{maximum_entropy}), we have
\begin{align*}
H(\w) - H(\vct{\pi}) =  H(\vct{\pi}) - H(\vct{p}),
\end{align*}
with $\vct{p} \in \Delta^{k_{\infty}-1}$ with $p_l = \sum_{r \in V_l}\pi_{r}$ for all $1 \le l \le k_{\infty}$. 
We show that the right-hand side of this equality is strictly positive:
\begin{align*}
H(\vct{\pi}) - H(\vct{p}) &= - \sum_{r=1}^d \pi_r \log \pi_r + \sum_{l=1}^{k_{\infty}} \left(\sum_{r \in V_l}\pi_{r}\right) \log \left(\sum_{r \in V_l}\pi_{r}\right)\\
&= - \sum_{l=1}^{k_{\infty}} \sum_{r \in V_l } \pi_r \log\left(\frac{\pi_r}{p_l}\right)\\
&= - \sum_{l=1}^{k_{\infty}} p_l \sum_{r \in V_l } \frac{\pi_r}{p_l} \log\left(\frac{\pi_r}{p_l}\right) \\
&\ge - \sum_{l=1}^{k_{\infty}} p_l \log\left(\frac{\sum_{r \in V_l }\pi_r^2}{p_l^2}\right),\\
&\ge 0.
\end{align*}
We used Jensen's inequality on the concave function $x\mapsto \log x$, and the fact that $\sum_{r \in V_l} \pi_r^2 \le p_l\sum_{r \in V_l} \pi_r = p_l^2$ for all $l$. Moreover, since all coordinates of $\vct{\pi}$ are strictly positive, equality holds if and only if $\pi_r = p_l$ for all $l$ and $ r \in V_l$ which implies that the partition must be trivial; i.e., $k_\infty = d$. Recall that this does not happen since $E_\infty$ is non-empty. 

On the other hand, by setting $\epsilon$ sufficiently small (smaller than all the limits in the definition of $E_\infty$), any edge in $E_{\infty}$ will eventually (and permanently from then on) be in $E_n$. Therefore the number of connected components of $G_n$ does not exceed that of $G_\infty$:\ $k_n \le k_{\infty}$ for $n$ sufficiently large. We conclude that $\psi_n(\mtx{v}^{(n)}/n) > \psi_n(\matoverlap^{(n)}/n)$ for all $n$ large enough. Therefore $\matoverlap^{(n)}$ is not always a maximizer of $\psi_n$, and this leads to a contradiction.   
\end{proofof}

\section{Proof of Theorem~\ref{collision_asymptotic}}

Our proof is based on the method of Laplace from asymptotic analysis: when the entries of $\matoverlap$ are large, the sum defining $q(\matoverlap)$ is dominated by its largest term corrected by a sub-exponential term which is represented by a Gaussian integral (see, e.g., \cite{debruijn1970asymptotic} for the univariate case). Since we are in a  multivariate situation, the asymptotics of $q$ depend on which subset of the entries of $\matoverlap$ are large. Our approach is inspired by~\cite{achlioptas2004two_journal}.        
 We recall that for $\matoverlap \in \Zp^{d \times d}$, 
\[q(\matoverlap) = \sum_{\stackrel{\mateulerflow \in \Zp^{d \times d} \cap \Fspace}{0 \le \eulerflow_{rs} \le \overlap_{rs}}} \left( \prod_{r,s=1}^d {{\overlap_{rs}}\choose{\eulerflow_{rs}}} \alpha^{\eulerflow_{rs}} (1-\alpha)^{\overlap_{rs} - \eulerflow_{rs}}\right). \]
Let $G = (V,E)$ with $V=\{1,\cdots,d\}$ and $E = \{(r,s) \in V^2~:~ r\neq s\}$. The graph $G$ will be used to store information about which entries of $\matoverlap$ are going to infinity linearly in $n$, and which entries are not. We can split the sum defining $q$ into a double sum, one involving the large terms ($A$ in subsequent notation), and the rest:
\[q(\matoverlap) = \sum_{\substack{0\le \eulerflow_{rs}' \le \overlap_{rs} \\ (r,s) \notin E}} \prod_{(r,s) \notin E} {{\overlap_{rs}}\choose{\eulerflow_{rs}'}} \alpha^{\eulerflow_{rs}'} (1-\alpha)^{\overlap_{rs} - \eulerflow_{rs}'} A(\mateulerflow',\matoverlap),\]
with
\[A(\mateulerflow',\matoverlap) = \sum_{\substack{0\le \eulerflow_{rs} \le \overlap_{rs} \\ (r,s) \in E}} \prod_{(r,s) \in E} {{\overlap_{rs}}\choose{\eulerflow_{rs}}} \alpha^{\eulerflow_{rs}} (1-\alpha)^{\overlap_{rs} - \eulerflow_{rs}}
\indi\left\{\bm{M}_G(\mateulerflow,\mateulerflow') \in \Fspace\right\}, \]
where for two $d \times d$ matrices $\mtx{a}, \mtx{b}$, $\bm{M}_G(\mtx{a} , \mtx{b})$ is the $d \times d$ matrix with entries $a_{rs}$ if $(r,s)\in E$ and $b_{rs}$ otherwise. The quantity $A$ will be approximated using the Laplace method. Recall from the expressions~\eqref{varphi_function} and~\eqref{varphi_minimization}  that 
\[\varphi_{\matoverlap} (\x) = \sum_{(r,s)\in E} \overlap_{rs} \kull{x_{rs}}{\alpha},\] 
and  
\[\vartheta(\mateulerflow,\matoverlap) =  \underset{\substack{\x \in [0,1]^{d \times d}\\ \bm{M}_G(\x\odot \matoverlap , \mateulerflow) \in \Fspace}}{\min} ~\varphi_{\matoverlap}(\x).\]
Let $\x^*(\mateulerflow,\matoverlap)$ be the optimal solution of the above optimization problem.   

Before stating our asymptotic approximation result for $A$, we state an important lemma on the boundedness of the entries of $\x^*(\mateulerflow,\matoverlap)$, where the bounds depend only on $\epsilon$ and $\alpha$.

\begin{lemma}\label{x-boundedness}
Let $G$ be fixed as above, $\alpha \in (0,1)$ and $\epsilon \in (0,1)$. There exist two constants $0 <c_l \le c_u < 1$ depending only on $d$, $\alpha$ and $\epsilon$ such that the following is true: For all integers $n \ge 1$, and $\matoverlap \in \{0,\cdots,n\}^{d \times d}$ such that $\overlap_{rs} \ge \epsilon n$ iff $(r,s) \in E$. For all $\mateulerflow'  \in \{0,\cdots,n\}^{\bar{E}}$ such that $0 \le \eulerflow'_{rs} \le \overlap_{rs}$ for all $(r,s) \notin E$, we have 
\[ c_l \le \min_{(r,s)\in E} x^*_{rs} \le \max_{(r,s)\in E}x^*_{rs} \le c_u.\]
\end{lemma}
Therefore, the entries of $\x^*$ can effectively be treated as constants throughout the rest of the proof. Now we state our asymptotic estimate for $A$.

\begin{proposition}\label{laplace}
 Let $G$ be fixed as above, and $\epsilon>0$. For all $n$ sufficiently large, all $\matoverlap \in \{0,\cdots,n\}^{d \times d}$ with $\overlap_{rs} \ge \epsilon n$ iff $(r,s) \in E$, and all $\mateulerflow'  \in \{0,\cdots,n\}^{\bar{E}}$ such that $0 \le \eulerflow'_{rs} \le \overlap_{rs}$ for all $(r,s) \notin E$, we have
\begin{align*}
 A(\mateulerflow',\matoverlap) ~~\asymp_{G,d,\epsilon,\alpha} \frac{e^{-\vartheta(\mateulerflow',\matoverlap)}}{P_{G}(\matoverlap)^{1/2}}. 
\end{align*}
Here, the symbol $``\asymp_{G,d,\epsilon,\alpha}"$ means that the ratio is upper- and lower-bounded by constants depending only on $G$, $d$, $\epsilon$ and $\alpha$. 
\end{proposition} 

By the above proposition, we have
\[q(\matoverlap) ~~\asymp_{G,d,\epsilon,\alpha} \sum_{\stackrel{\mateulerflow \in \Zp^{\bar{E}}}{0 \le \eulerflow_{rs} \le \overlap_{rs}}} \left( \prod_{(r,s)\notin E} {{\overlap_{rs}}\choose{\eulerflow_{rs}}} \alpha^{\eulerflow_{rs}} (1-\alpha)^{\overlap_{rs} - \eulerflow_{rs}}\right) \frac{e^{-\vartheta(\mateulerflow,\matoverlap)}}{P_{G}(\matoverlap)^{1/2}}.\]
The estimate above (ignoring the term $P_G(\matoverlap)$) can be interpreted as the expected value of the function $e^{-\vartheta(\mateulerflow,\matoverlap)}$ under the law of the random variable $\mateulerflow$ where each entry $\eulerflow_{rs}$ for $(r,s) \notin E$ is independently binomial with parameters $\alpha$ and $\overlap_{rs}$. From here, the bounds claimed in Theorem~\ref{collision_asymptotic} follow immediately.

\vspace{.3cm}
\noindent \begin{proofof}{Proposition~\ref{laplace}}
 We will show that 
\[ A(\mateulerflow',\matoverlap) ~~\asymp_{G,d,\epsilon,\alpha} ~~e^{-\vartheta(\mateulerflow',\matoverlap)}~ \prod_{(r,s)\in E}\overlap_{rs}^{-1/2}~\int_{\Fspace(G)} e^{- \sum_{(r,s)\in E} z_{rs}^2/2\overlap_{rs}} ~\mathrm{d}\z.\]
Then the result follows by applying Proposition~\ref{gaussian-integral} to evaluate the Gaussian integral. We proceed by showing the upper and lower bounds separately.

\paragraph{The upper bound} 
For $\mateulerflow' \in \Zp^{\bar{E}}, \matoverlap \in \Zp^{d\times d}$ fixed and some parameter $C(\matoverlap)>0$ to be adjusted, let 
\[\Omega = \left\{ \mateulerflow \in \Z_+^{E} ~:~ \bm{M}_G(\mateulerflow,\mateulerflow') \in \Fspace,~ 0 \le \eulerflow_{rs} \le \overlap_{rs},~ \sum_{(r,s) \in E} \frac{(\eulerflow_{rs} - x^*_{rs}\overlap_{rs})^2}{x^*_{rs}(1-x^*_{rs})\overlap_{rs}} \le C(\matoverlap)^2\right\}.\] 

For $\mateulerflow \in \Zp^{E}$, we let $\x \in [0,1]^{E}$ defined by $x_{rs} = \eulerflow_{rs}/\overlap_{rs}$ for all $(r,s)\in E$. 
We upper bound the binomial coefficients ${{\overlap_{rs}}\choose{\eulerflow_{rs}}}$ based on whether $\mateulerflow$ is in $\Omega$ or not:
\begin{itemize}
\item If $\mateulerflow \in \Omega$ we use the upper bound ${{\overlap_{rs}}\choose{\eulerflow_{rs}}} \le \left(2\pi \overlap_{rs} x_{rs}(1-x_{rs})\right)^{-1/2}\exp \overlap_{rs} H(x_{rs})$. 
\item Otherwise, we use the upper bound ${{\overlap_{rs}}\choose{\eulerflow_{rs}}} \le 3\sqrt{\overlap_{rs}} \exp \overlap_{rs} H(x_{rs})$.
\end{itemize}
Here, $H(x_{rs}) = -x_{rs}\log x_{rs} - (1-x_{rs})\log(1-x_{rs})$. Thus, the summand in $A(\mateulerflow',\matoverlap)$ is bounded by 
\begin{align*}
\prod_{(r,s)\in E}& \left(2\pi \overlap_{rs} x_{rs}(1-x_{rs})\right)^{-1/2} \exp \overlap_{rs} \kull{x_{rs}}{\alpha} 
= \prod_{(r,s) \in E} \left(2\pi \overlap_{rs} x_{rs}(1-x_{rs})\right)^{-1/2} \exp (-\varphi_{\matoverlap}(\x))
\end{align*}
if $\mateulerflow \in \Omega$, and \[\prod_{(r,s)\in E} 3\overlap_{rs}^{1/2} \exp (-\varphi_{\matoverlap}(\x))\] if not.
The function $\varphi_{\matoverlap}$ is smooth, and we have $\frac{\mathrm{d}\varphi_{\matoverlap}}{\mathrm{d}x_{rs}}(\x) = \overlap_{rs}\log\left(\frac{x_{rs}(1-\probquery)}{\probquery (1-x_{rs})}\right)$, and $\frac{\mathrm{d}^2\varphi_{\matoverlap}}{\mathrm{d}x_{rs}^2}(\x) = \frac{\overlap_{rs}}{x_{rs}(1-x_{rs})} \ge 0$. Therefore, by convexity,
\[\varphi_{\matoverlap}(\x)  \ge \varphi_{\matoverlap}(\x^*) + \frac{1}{2}\sum_{(r,s)\in E}^d \frac{\overlap_{rs}}{x^*_{rs}(1-x^*_{rs})} (x_{rs} - x^*_{rs})^2.\]    
Let 
\begin{equation} \label{quadratic}
\ell_{\matoverlap}(\mateulerflow) = \sum_{(r,s)\in E} \frac{(\eulerflow_{rs} - \overlap_{rs}x^*_{rs})^2}{x^*_{rs}(1-x^*_{rs})\overlap_{rs}}.
\end{equation}
Based on Lemma~\ref{x-boundedness}, all the entries of $\x^*$ will be treated as constants. If $\mateulerflow \in \Omega$ then $\frac{(\eulerflow_{rs} - \overlap_{rs}x^*_{rs})^2}{x^*_{rs}(1-x^*_{rs})\overlap_{rs}} \le \ell_{\matoverlap}(\mateulerflow) \le C(\matoverlap)^2$, and  
\[\eulerflow_{rs} \in \left[\overlap_{rs}x^*_{rs} \pm C(\matoverlap) \sqrt{x^*_{rs}(1-x^*_{rs})\overlap_{rs}}\right].\] 
Now let us assume that $C(\matoverlap) = o(\overlap_{rs}^{1/2})$ for all $(r,s)\in E$. Then we have
\[\frac{\eulerflow_{rs}}{\overlap_{rs}}(1-\frac{\eulerflow_{rs}}{\overlap_{rs}}) \ge x_{rs}^*(1-x_{rs}^*) - o_{n}(1).\] 
If $\mateulerflow \notin \Omega$ then $\ell(\mateulerflow) \ge C(\matoverlap)^2$, therefore $A$ is bounded by 
\begin{equation} \label{large_sum} 
\left(\prod_{(r,s)\in E} (2\pi \overlap_{rs} x^*_{rs}(1-x^*_{rs}))^{-1/2} \sum_{\mateulerflow \in \Omega} e^{-\ell_{\matoverlap}(\mateulerflow)/2}
+
\prod_{(r,s)\in E} 3\overlap_{rs}^{1/2} \sum_{\substack{\mateulerflow \notin \Omega \\ 0 \le \eulerflow_{rs} \le \overlap_{rs}}}  e^{-C(\matoverlap)^2/2} \right)
\cdot \exp(-\varphi_{\matoverlap}(\x^*)).
\end{equation} 
The second term in the sum above is bounded by $c^{d^2}(\prod_{(r,s)\in E} \overlap_{rs}^{3/2}) e^{-C(\matoverlap)^2/2}$ for some constant $c>0$. Taking $C(\matoverlap)^2 = 5 \log \prod_{(r,s)\in E} \overlap_{rs}$, this term is $c^{d^2}\prod_{(r,s)\in E} \overlap_{rs}^{-1} = \bigo(n^{-|E|})$. Moreover, for all $(r,s) \in E$, $\overlap_{rs} > \epsilon n$, therefore 
\[C(\matoverlap)^2 \le 4d^2\log n \ll  \overlap_{rs};\] 
this choice satisfies the condition $C(\matoverlap) = o(\overlap_{rs}^{1/2})$ for $(r,s) \in E$.
On the other hand, let 
\begin{equation}\label{gaussian_sum}
S(\matoverlap) = \sum_{\mateulerflow \in \Z^{E}} \indi\{\bm{M}_G(\mateulerflow,\mateulerflow') \in \Fspace\}\exp(-\ell_{\matoverlap}(\mateulerflow)/2),
\end{equation}
with $\ell_{\matoverlap}$ defined in~\eqref{quadratic}. We ignored the dependence of $S$ on $\mateulerflow'$ in the notation on purpose: this dependence is inessential.   
The first sum in~\eqref{large_sum} is upper bounded by $S(\matoverlap)$. 
Therefore 
\begin{equation}\label{collision_upper_bound}
A(\mateulerflow',\matoverlap) \le c_u \prod_{(r,s)\in E} \overlap_{rs}^{-1/2}~(S(\matoverlap) +\varepsilon_n)~e^{-\vartheta(\mateulerflow',\matoverlap)}, 
\end{equation}
for some $c_u$ depending on $\epsilon$ and $d$, and  $\epsilon_n \to 0 $ as $n \to \infty$.
We now turn our attention to the lower bound, deferring the analysis of the Gaussian sum $S(\matoverlap)$ to a subsequent paragraph.   

\paragraph{The lower bound}
We have 
\[A(\mateulerflow',\matoverlap) \ge \sum_{\mateulerflow \in \Omega} \left( \prod_{(r,s)\in E} {{\overlap_{rs}}\choose{\eulerflow_{rs}}} \probquery^{\eulerflow_{rs}} (1-\probquery)^{\overlap_{rs} - \eulerflow_{rs}}\right).\] 
Using ${{\overlap_{rs}}\choose{\eulerflow_{rs}}} \ge (8\pi \eulerflow_{rs}(1-\eulerflow_{rs}/\overlap_{rs}))^{-1/2}e^{H(\eulerflow_{rs}/\overlap_{rs})}$ for all $(r,s)\in E$, we get 
\[A(\mateulerflow',\matoverlap) \ge \sum_{\mateulerflow \in \Omega}  \left(\prod_{(r,s)\in E}8\pi \eulerflow_{rs}(1-\frac{\eulerflow_{rs}}{\overlap_{rs}})\right)^{-1/2} \cdot \exp\left(-\varphi_{\matoverlap}\left(\mateulerflow/\matoverlap\right)\right).\] 
For $\mateulerflow \in \Omega$, we have $\frac{\eulerflow_{rs}}{\overlap_{rs}}(1-\frac{\eulerflow_{rs}}{\overlap_{rs}}) \le x_{rs}^*(1-x_{rs}^*) + o_{n}(1)$ for all $r,s$, and since $\varphi_{\matoverlap}$ is a smooth function, a Taylor expansion yields 
\[\varphi_{\matoverlap}(\mateulerflow/\matoverlap) = \varphi_{\matoverlap}(\x^*) + \frac{1}{2}\sum_{(r,s)\in E} \frac{(\eulerflow_{rs} - x^*_{rs}\overlap_{rs})^2}{x^*_{rs}(1-x^*_{rs})\overlap_{rs}} + o_n(1).\]
Therefore,
\begin{align*}
A(\mateulerflow',\matoverlap) &\ge c_l e^{-\varphi_{\matoverlap}(\x^*)}\prod_{(r,s)\in E} \overlap_{r,s}^{-1/2} \cdot \sum_{\mateulerflow \in \Omega} \exp\left( - \frac{1}{2}\sum_{rs} \frac{(\eulerflow_{rs} - x^*_{rs}\overlap_{rs})^2}{x^*_{rs}(1-x^*_{rs})\overlap_{rs}}\right)\\
&= c_l e^{-\vartheta(\mateulerflow',\matoverlap)}\prod_{(r,s)\in E} \overlap_{r,s}^{-1/2} \cdot \left(S(\matoverlap) - \varepsilon_n\right),
\end{align*}
where $c_l = c_l(\epsilon,d)$, $S(\matoverlap)$ is defined in~\eqref{gaussian_sum} and 
\begin{align*}
\varepsilon_n := \prod_{(r,s)\in E}(\overlap_{rs}+1)e^{-C(\matoverlap)^2/2} + \sum_{\substack{\mateulerflow \in \Zp^{E} \\ \eulerflow_{rs} \ge \overlap_{rs}+1}} \exp - \ell_{\matoverlap}(\mateulerflow)/2 
+ \sum_{\mateulerflow \in \Z_{-}^{E}} \exp - \ell_{\matoverlap}(\mateulerflow)/2.
\end{align*}
We take $C(\matoverlap)^2 = 4 \log \prod_{(r,s)\in E} \overlap_{rs}$. This makes the first term in $\varepsilon_n$ bounded by $\prod_{r,s} \overlap_{rs}^{-1} = \bigo (n^{-|E|})$. On the other hand, the remaining tail sums are easily bounded by the tail probability function of a normal random variable (i.e., the error function):
\begin{align*}
\sum_{\substack{\mateulerflow \in \Zp^{E} \\ \eulerflow_{rs} \ge \overlap_{rs}+1}} \exp - \ell_{\matoverlap}(\mateulerflow)/2 &\le 
\prod_{(r,s)\in E} \overlap_{rs}^{1/2} \erfc\left( \sqrt{\frac{x^*_{rs}}{1-x^*_{rs}}\overlap_{rs}}\right),\\
\sum_{\mateulerflow \in \Z_{-}^{E}} \exp - \ell_{\matoverlap}(\mateulerflow)/2 &\le 
\prod_{(r,s)\in E} \overlap_{rs}^{1/2}  \erfc\left( \sqrt{\frac{1-x^*_{rs}}{x^*_{rs}}\overlap_{rs}} - \frac{1}{\sqrt{x^*_{rs}(1-x^*_{rs})\overlap_{rs}}}\right),
\end{align*}
with $\erfc(x) = \int_{x}^\infty e^{-t^2/2}dt$. Since $\erfc(x) \le e^{-x^2/2}/x$ for all $x >0$, these two terms decay in a sub-Gaussian way in $n$.

\paragraph{Bounding the Gaussian sum.} Here we approximate $S$ by a continuous Gaussian integral. We prove that 
\begin{align*}
S(\matoverlap) \asymp \int_{\Fspace(G)} \exp\left(- \sum_{(r,s)\in E} \frac{\overlap_{rs}^{-1}}{2x^*_{rs}(1-x^*_{rs})}z_{rs}^2\right) \mathrm{d}\z, 
\end{align*}
where the symbol $``\asymp"$ hides constants depending on $G,\epsilon,d$ and $\alpha$ as $n \to \infty$. 
For $\mateulerflow \in \Fspace(G)$ an array of integer numbers such that $0 \le \eulerflow_{rs} \le \overlap_{rs}$, let $T(\mateulerflow) = \mateulerflow + \mathcal{C} \cap \Fspace(G)$ where $\mathcal{C} = [-1/2,1/2]^{E}$. The sum is understood in the Minkowski sense. $T(\mateulerflow)$ is a ``tile" of side 1 centered around $\mateulerflow$. Two crucial facts are $(i)$ : $T(\mateulerflow)$ and $T(\mateulerflow')$ are of disjoint interiors when $\mateulerflow \neq \mateulerflow'$ and $(ii)$ : $T(\mateulerflow) \subset \Fspace(G)$. Now for a fixed $\mateulerflow$, let $\z \in T(\mateulerflow)$. For $r,s \in E$, we have $\eulerflow_{rs}-1/2 \le z_{rs} \le \eulerflow_{rs}+1/2$ and $\frac{\eulerflow_{rs}-1/2}{\overlap_{rs}} -x^*_{rs} \le \frac{z_{rs}}{\overlap_{rs}} -x^*_{rs} \le \frac{\eulerflow_{rs}+1/2}{\overlap_{rs}} -x^*_{rs}$. Thus
 \[\left(\frac{z_{rs}}{\overlap_{rs}} -x^*_{rs}\right)^2 \le \max\left\{\left(\frac{\eulerflow_{rs}-1/2}{\overlap_{rs}} -x^*_{rs}\right)^2,\left(\frac{\eulerflow_{rs}+1/2}{\overlap_{rs}} -x^*_{rs}\right)^2\right\}.\]
Using the fact $\max\{(x-1/2)^2,(x+1/2)^2\} \le 2x^2 +1$ for all $x \in \R$, we get
 \[\exp\left(- \sum_{(r,s)\in E} \frac{\overlap_{rs}^{-1}}{4x^*_{rs}(1-x^*_{rs})} \left(2\left(\eulerflow_{rs} - \overlap_{rs}x^*_{rs}\right)^2+1\right)\right) \le \exp\left(- \sum_{(r,s)\in E} \frac{\overlap_{rs}^{-1}}{4x^*_{rs}(1-x^*_{rs})} \left(z_{rs} - \overlap_{rs}x^*_{rs}\right)^2\right).\]

By integrating both sides of the above inequality on $T(\mateulerflow)$ in the variable $\z$, and summing over all $\mateulerflow$ with integer entries such that $\bm{M}_{G}(\mateulerflow,\mateulerflow') \in \Fspace$, we get 
\begin{align*}
\vol(\mathcal{C}\cap \Fspace(G))e^{-\sum_{(r,s)\in E} \frac{\overlap_{rs}^{-1}}{4x^*_{rs}(1-x^*_{rs})}} ~S(\matoverlap) 
&\le \sum_{\substack{\mateulerflow \in \Z^{E}\\\bm{M}_{G}(\mateulerflow,\mateulerflow') \in \Fspace}}\int_{T(\mateulerflow)} \exp\left(- \sum_{(r,s)\in E} \frac{\overlap_{rs}^{-1}\left(z_{rs} - \overlap_{rs}x^*_{rs}\right)^2}{4x^*_{rs}(1-x^*_{rs})} \right) \mathrm{d}\z,\\
&= \sum_{\substack{\mateulerflow \in \Z^{E}\\\bm{M}_{G}(\mateulerflow,\mateulerflow') \in \Fspace}} \int_{T(\mateulerflow-\x^*\odot\matoverlap)} \exp\left(- \sum_{(r,s)\in E} \frac{\overlap_{rs}^{-1}}{4x^*_{rs}(1-x^*_{rs})} z_{rs}^2\right) \mathrm{d}\z,
\end{align*}
where $\vol$ is the volume according to the $\dim(\Fspace(G))$-dimensional Lebesgue measure. Since $\bm{M}_{G}(\x^*\odot \matoverlap,\mateulerflow') \in \Fspace$, we have $\mateulerflow - \x^*\odot \matoverlap \in \Fspace$ for all $\mateulerflow$ we are summing over. Moreover, since the tiles $T(\mateulerflow)$ are of mutually disjoint interiors, and given that their union is in $\Fspace(G)$, the left-hand side is upper bounded by (there is actually equality)
\[\int_{\Fspace(G)} \exp\left(- \sum_{(r,s)\in E} \frac{\overlap_{rs}^{-1}}{4x^*_{rs}(1-x^*_{rs})} z_{rs}^2\right) \mathrm{d}\z. \] 

Here, to get sharper constants, one could apply a theorem by Vaaler~\cite{vaaler1979} which states that the volume of any linear subspace intersected with the cube $\mathcal{C}$ is at least 1; i.e., $\vol(\mathcal{C}\cap \Fspace(G)) \ge 1$. This yields
\[S(\matoverlap) \le e^{\sum_{(r,s)\in E} \frac{\overlap_{rs}^{-1}}{4x^*_{rs}(1-x^*_{rs})}} \cdot \int_{\Fspace(G)} \exp\left(- \sum_{(r,s)\in E} \frac{\overlap_{rs}^{-1}}{4x^*_{rs}(1-x^*_{rs})}z_{rs}^2\right) \mathrm{d}\z.\]   

As for the reverse inequality, slightly more care is needed in constructing the approximation. For a given $\mateulerflow$, let $\Omega^+ = \{(r,s) ~:~ \eulerflow_{rs} \ge x^*_{rs}\overlap_{rs}+1/2\}$ and $\Omega^- = \{(r,s) ~:~ \eulerflow_{rs} \le x^*_{rs}\overlap_{rs}-1/2\}$. For $\z \in T(\mateulerflow)$, 
we have $(z_{rs} - x^*_{rs}\overlap_{rs})^2 \ge (\eulerflow_{rs}-x^*_{rs}\overlap_{rs}-1/2)^2$ if $(r,s)\in \Omega^{+}$ and $(z_{rs} - x^*_{rs}\overlap_{rs})^2 \ge (\eulerflow_{rs}-x^*_{rs}\overlap_{rs} + 1/2)^2$ if $(r,s)\in \Omega^{-}$. Otherwise, for $(r,s) \notin \Omega^+ \cup \Omega^-$, we have $\abs{\eulerflow_{rs}- x^*_{rs}\overlap_{rs}} < 1/2$ and $\abs{(z_{rs} - x^*_{rs}\overlap_{rs})^2 - (\eulerflow_{rs} - x^*_{rs}\overlap_{rs})^2} < 1/2(1+1/2) = 3/4$. Therefore
\begin{align*}
&\sum_{(r,s) \in \Omega^+} \frac{(\eulerflow_{rs} - x^*_{rs}\overlap_{rs}+1/2)^2}{\overlap_{rs}x^*_{rs}(1-x^*_{rs})} + \sum_{(r,s) \in \Omega^-} \frac{(\eulerflow_{rs} - x^*_{rs}\overlap_{rs}-1/2)^2}{\overlap_{rs}x^*_{rs}(1-x^*_{rs})} + \sum_{(r,s) \notin \Omega^+\cup \Omega^-} \frac{(\eulerflow_{rs} - x^*_{rs}\overlap_{rs})^2}{\overlap_{rs}x^*_{rs}(1-x^*_{rs})}\\
&\le \sum_{(r,s)\in E} \frac{(z_{rs} - x^*_{rs}\overlap_{rs})^2}{\overlap_{rs}x^*_{rs}(1-x^*_{rs})} + \sum_{(r,s)\in E}\frac{3\overlap_{rs}^{-1}}{4x^*_{rs}(1-x^*_{rs})}.
\end{align*}
On the other hand, $(\eulerflow_{rs} - x^*_{rs}\overlap_{rs})^2 \le (\eulerflow_{rs} - x^*_{rs}\overlap_{rs}+1/2)^2$ when $(r,s)\in \Omega^+$ and $(\eulerflow_{rs} - x^*_{rs}\overlap_{rs})^2 \le (\eulerflow_{rs} - x^*_{rs}\overlap_{rs}-1/2)^2$ when $(r,s)\in \Omega^-$.
After integrating on $T(\mateulerflow)$ and summing over all $\mateulerflow \in \Z^{E}$ such that $\mateulerflow - \x^* \odot \matoverlap \in \Fspace$, we obtain:
\[\vol(\mathcal{C} \cap \Fspace(G)) S(\matoverlap) \ge e^{-\sum_{(r,s)\in E}\frac{3\overlap_{rs}^{-1}}{8x^*_{rs}(1-x^*_{rs})}} \cdot \sum_{\substack{\mateulerflow \in \Z^{E}\\ \mateulerflow - \x^* \odot \matoverlap \in \Fspace}} \int_{T(\mateulerflow - \x^* \odot \matoverlap)} \exp\left(- \sum_{(r,s)\in E} \frac{\overlap_{rs}^{-1}}{2x^*_{rs}(1-x^*_{rs})}z_{rs}^2\right) \mathrm{d}\z,\]
and the last sum is equal to 
\[\sum_{\mateulerflow \in (\Z^{E}+\x^* \odot \matoverlap)\cap \Fspace(G)} \int_{T(\mateulerflow)} \exp\left(- \sum_{(r,s)\in E} \frac{\overlap_{rs}^{-1}}{2x^*_{rs}(1-x^*_{rs})}z_{rs}^2\right) \mathrm{d}\z = \int_{\Fspace(G)} \exp\left(- \sum_{(r,s)\in E} \frac{\overlap_{rs}^{-1}}{2x^*_{rs}(1-x^*_{rs})}z_{rs}^2\right) \mathrm{d}\z.\] 
Finally,
\[ S(\matoverlap) \ge c(G,d)e^{- \sum_{(r,s)\in E} \frac{3\overlap_{rs}^{-1}}{2x^*_{rs}(1-x^*_{rs})}} \cdot \int_{\Fspace(G)}\exp\left(- \sum_{(r,s)\in E} \frac{\overlap_{rs}^{-1}}{2x^*_{rs}(1-x^*_{rs})}z_{rs}^2\right) \mathrm{d}\z.\]
\end{proofof}

\vspace{.5cm}
\begin{proofof}{Lemma~\ref{x-boundedness}}
Recall that $\x^*$ is the unique minimizer of the function \[\varphi_{\matoverlap} = \sum_{(r,s)\in E} \overlap_{rs}\kull{x_{rs}}{\alpha}\] 
on $[0,1]^{d\times d}$ subject to $\bm{M}_{G}(\x^*\odot \matoverlap,\mateulerflow) \in \Fspace$.  Recall also that the entries of $\x^*$ admit the expressions
\[x^*_{rs} = \frac{\alpha}{\alpha +(1-\alpha)e^{\lambda_{r}^* - \lambda_{s}^*}},\]  
for all $ (r,s) \in E$. The vector $\vct{\lambda}^* \in \R^d$ is the unique solution up to global shifts to the dual optimization problem (strong duality holds here~\cite{boyd02,rockafellar70}) 
\begin{align}\label{dual-problem}
\sup_{\vct{\lambda} \in \R^d} \left\{ \sum_{(r,s)\notin E} \eulerflow_{rs}(\lambda_{r} - \lambda_{s}) + \sum_{(r,s)\in E} \overlap_{rs}\log\left(\frac{e^{\lambda_{r} - \lambda_{s}}}{\alpha + (1-\alpha)e^{\lambda_{r} - \lambda_{s}}}\right) \right\}.
\end{align}
Our claim reduces to the boundedness of the differences $|\lambda_r^* - \lambda_s^*|$ for all $(r,s)\in E$ independently of $n,\matoverlap,\mateulerflow$ and $r,s$. We will shortly prove the following inequality
\begin{align}\label{lambda-boundedness}
\sum_{(r,s)\in E} \overlap_{rs} (\lambda_r^* - \lambda_s^*)^2 \le \kappa(\alpha) \sum_{(r,s)\in E} \overlap_{rs},
\end{align}
where $\kappa(\alpha) = \frac{1}{\alpha^2} + \frac{1}{(1-\alpha)^2}$. Assuming the above is true, by the Cauchy-Schwarz inequality, we would have
\begin{align*}
\sum_{(r,s)\in E} |\lambda_r^* - \lambda_s^*| \le \left(\sum_{(r,s)\in E} \overlap_{rs}^{-1} \right)^{1/2}\left(\kappa(\alpha) \sum_{(r,s)\in E} \overlap_{rs} \right)^{1/2}
\le d^2(\kappa(\alpha)/\epsilon)^{1/2},
\end{align*}
since $\epsilon n \le \overlap_{rs} \le n$ for all $(r,s)\in E$. We would then be done. 
Now, the inequality~\eqref{lambda-boundedness} follows from convexity considerations. We let $\phi$ be the function being maximized in~\eqref{dual-problem}. By concavity of $\phi$, we have
\begin{align}\label{concavity}
\phi(\vct{\lambda}^*) - \phi(\vct{0}) \le \vct{\lambda}^{*\intercal}\nabla\phi(\vct{0}) + \frac{1}{2} \vct{\lambda}^{*\intercal}\nabla^2\phi(\vct{0})\vct{\lambda}^*.
\end{align}
The gradient and the Hessian of $\phi$ are 
\[ \left[\nabla\phi(\vct{\lambda})\right]_r = \sum_{s: (r,s)\in E} \frac{\alpha\overlap_{rs}}{\alpha+(1-\alpha)e^{\lambda_r-\lambda_s}} -\frac{\alpha\overlap_{sr}}{\alpha+(1-\alpha)e^{\lambda_s-\lambda_r}} + \sum_{s: (r,s)\notin E} \eulerflow_{rs} - \eulerflow_{sr}, \quad r \in \{1,\cdots,d\},\]
\[\nabla^2\phi(\vct{\lambda}) = -\alpha(1-\alpha)\sum_{(r,s)\in E} w_{rs}(\vct{\lambda})(\vct{e}_r - \vct{e}_s)(\vct{e}_r-\vct{e}_s)^\intercal,\]
with \[w_{rs}(\vct{\lambda}) = \frac{\overlap_{rs}e^{\lambda_r-\lambda_s}}{(\alpha+(1-\alpha)e^{\lambda_r-\lambda_s})^2} + \frac{\overlap_{sr}e^{\lambda_s-\lambda_r}}{(\alpha+(1-\alpha)e^{\lambda_s-\lambda_r})^2},\]
and $\vct{e}_1,\cdots,\vct{e}_d$ being the standard unit vectors in $\R^d$.
The concavity inequality~\eqref{concavity} becomes
\[ \phi(\vct{\lambda}^*) \le \alpha \sum_{(r,s)\in E} \overlap_{rs}(\lambda_r^*  - \lambda_s^*) + \sum_{(r,s)\notin E} \eulerflow_{rs}(\lambda_r^*  - \lambda_s^*)- \alpha(1-\alpha) \sum_{(r,s)\in E}\overlap_{rs}(\lambda_r^*  - \lambda_s^*)^2.\]
Substituting in the expression of $\phi(\vct{\lambda}^*)$, the term $\sum_{(r,s)\notin E} \eulerflow_{rs}(\lambda_r^*  - \lambda_s^*)$ cancels out on both sides and we get
\[
\sum_{(r,s)\in E} \overlap_{rs}\log\left(\frac{e^{\lambda_{r}^* - \lambda_{s}^*}}{\alpha + (1-\alpha)e^{\lambda_{r}^* - \lambda_{s}^*}}\right)
\le 
\alpha \sum_{(r,s)\in E} \overlap_{rs}(\lambda_r^*  - \lambda_s^*) - \alpha(1-\alpha) \sum_{(r,s)\in E}\overlap_{rs}(\lambda_r^*  - \lambda_s^*)^2,
\]
which can be written as 
\begin{align}\label{concavity_3}
 \sum_{(r,s)\in E}\overlap_{rs}\left(\alpha(1-\alpha)(\lambda_r^*  - \lambda_s^*)^2 +(1-\alpha)(\lambda_r^*  - \lambda_s^*) - \log\left(\alpha + (1-\alpha)e^{\lambda_{r}^* - \lambda_{s}^*}\right)\right) \le 0.
 \end{align}
Now we approximate the logarithm by the positive part: $\log(\alpha + (1-\alpha)e^{x}) \le x_+ = \max\{0,x\}$ for all $x \in \R$ and $\alpha \in (0,1)$, so that we almost get a quadratic polynomial inequality. We make this a genuine quadratic inequality by applying the additional approximation that for all $x\in \R$ and $\alpha \in (0,1)$: 
\[\alpha(1-\alpha)x^2 + (1-\alpha)x - x_+ \ge \frac{\alpha(1-\alpha)}{2}x^2 - \frac{1-\alpha}{2\alpha} - \frac{\alpha}{2(1-\alpha)}.\]
This is easy to check by verifying that the discriminants of the resulting quadratics (one for $x\ge 0$ and one for $x <0$) are negative.  
Now, inequality~\eqref{concavity_3} implies 
\[\frac{\alpha(1-\alpha)}{2}\sum_{(r,s) \in E} \overlap_{rs}(\lambda_r^*  - \lambda_s^*)^2 \le \left(\frac{1-\alpha}{2\alpha} + \frac{\alpha}{2(1-\alpha)}\right)\sum_{(r,s) \in E} \overlap_{rs}.\]
In other words,
\[\sum_{(r,s) \in E} \overlap_{rs}(\lambda_r^*  - \lambda_s^*)^2 \le \kappa(\alpha) \sum_{(r,s) \in E} \overlap_{rs}.\]
\end{proofof}

\section{Two proofs of Proposition~\ref{gaussian-integral}}
We first reduce the proof to the case where $G=K_d$ by a limiting argument. Let $G=(V,E)$ be a graph on $d$ vertices. If $G$ is not connected then  the constraints defining the space $\Fspace(G)$ decouple across the connected components of $G$ and so does the integrand $\exp -\half \sum_{(r,s)\in E} x_{rs}^2/w_{rs}$, therefore the Gaussian integral factors across the connected components of $G$. Hence, we may assume that $G$ is connected. Now, if 
\[\int_{\Fspace} ~e^{-\frac{1}{2}\sum_{rs} x_{rs}^2/w_{rs}} ~\textup{d}\x = (2\pi)^{((d-1)^2+d)/2} \left(\frac{\prod_{r,s} w_{rs}}{T(\w)}\right)^{1/2},\]
for all $\w \in \Rp^{d \times d}$ where $T = T_{K_{d}}$, then taking a limit $w_{rs} \to 0$ for all $(r,s) \notin E$, we get 
\[\frac{1}{\left(\prod_{(r,s)\notin E} w_{rs}\right)^{1/2}}\int_{\Fspace} ~e^{-\frac{1}{2}\sum_{rs} x_{rs}^2/w_{rs}} ~\mathrm{d}\x 
\longrightarrow 
c(G) \int_{\Fspace(G)} ~e^{-\frac{1}{2}\sum_{(r,s)\in E} x_{rs}^2/w_{rs}} ~\mathrm{d}\x,
\]
where $c(G)>0$ is a constant that only depends on $G$. On the other hand 
\[T(\w) \longrightarrow \frac{\nst(G)}{2^{d-1}d^{d-2}}~T_{G}(\w).\] 
Therefore  
\[c(G)\int_{\Fspace(G)} ~e^{-\frac{1}{2}\sum_{(r,s)\in E} x_{rs}^2/w_{rs}} ~\mathrm{d}\x =  (2\pi)^{((d-1)^2+d)/2} \left(\frac{2^{d-1}d^{d-2}}{\nst(G)} \frac{\prod_{(r,s) \in E} w_{rs}}{T_G(\w)}\right)^{1/2}.\]
Now we set $w_{rs} = 1$ for all $(r,s) \in E$ to clear out the constants. Since $\int_{\Fspace(G)} ~e^{-\frac{1}{2}\sum_{(r,s)\in E} x_{rs}^2} ~\mathrm{d}\x = (2\pi)^{\dim(\Fspace(G))/2}$, we get 
\[\int_{\Fspace(G)} ~e^{-\frac{1}{2}\sum_{(r,s)\in E} x_{rs}^2/w_{rs}} ~\mathrm{d}\x =  (2\pi)^{\dim(\Fspace(G))/2} \left(\frac{\prod_{(r,s) \in E} w_{rs}}{T_G(\w)}\right)^{1/2}.\]
Now it remains to prove the proposition for the complete graph.

\subsection{A combinatorial proof}
We proceed by adopting a combinatorial view on the structure of the space $\Fspace$. This will lead us to consider a very special basis of $\Fspace$ in which the computations become tractable. (Background on the concepts used in this construction can be found in~\cite{biggs1997algebraic}.) We first orient $K_d$ in such a way that every pair of distinct vertices is connected by two parallel edges pointing in opposite directions. There are $d(d-1)$ (oriented) edges in total. Then, the subgraphs whose edges are weighted by an array $\x \in \Fspace$ are called \emph{Eulerian}: the total weight of the incoming edges is equal to that of the outgoing edges on each vertex.
An important property of Eulerian graphs is that they can be decomposed into a superposition of cycles.  
In particular, fix a spanning tree $T^*$ of $K_d$ (the tree uses only one edge, if any, between each pair of vertices, and ignores their orientation). Every edge $e \notin T^*$ can be identified with the oriented cycle $C_e$ in the graph which consists of the oriented edge $e$ and the unique path between the endpoints of $e$ in the tree $T^*$ (where the direction of the edges on the path are flipped if necessary). Let $\vct{\chi}_e \in \{0,\pm1\}^{d(d-1)}$ be the indicator vector of the cycle $C_e$\footnote{Each non-zero entry in the vector corresponds to an edge present in the cycle, and the non-zero value is $+1$ if the cycle flows along the orientation of that edge, and $-1$ if the flow is in the opposite direction. In particular, the $e$th coordinate of $\vct{\chi}_e$ is always $+1$.}. Since a cycle is Eulerian, the vector $\vct{\chi}_e$---when folded into a $d \times d$ matrix---belongs to $\Fspace$. Furthermore, the family $\{\vct{\chi}_e ~:~ e \notin T^*\}$ is linearly independent since a cycle $C_e$ contains at least one edge---namely $e$---that is not contained in any other cycle $C_{e'}$, $e' \neq e$. There are exactly $d(d-1) - (d-1) = (d-1)^2$ off-tree edges in $K_d$, and this number coincides with the dimension of $\Fspace$. Therefore $\mathcal{B} = \{\vct{\chi}_e ~:~ e \notin T^*\}$ is a basis of $\Fspace$, that we henceforth call a \emph{fundamental cycle basis} of $\Fspace$. 

Let $\bm{P} \in \{0,\pm1\}^{(d-1)^2 \times d(d-1)}$ be the matrix where the rows are indexed by the off-tree edges of the graph, and whose $e$th row is equal to $\vct{\chi}_e$. The matrix $\bm{P}$ can be regarded as the \emph{cycle-edge incidence matrix} of the graph $K_d$: an entry $(e,e')$ is non-zero if and only if $e' \in C_e$.

Let $\bm{M} \in \R^{d(d-1) \times d(d-1)}$ be the diagonal matrix with entries $w_{rs}$, $r \neq s$ on the diagonal. Then by a change of variables
\begin{align*}
\int_{\Fspace} e^{-\sum_{rs} x_{rs}^2/2w_{rs}} \mathrm{d}\x &= Det(\bm{P}\bm{P}^\intercal)^{1/2} \int_{\R^{(d-1)^2}} e^{- \vct{z}^\intercal (\bm{P} \bm{M}^{-1} \bm{P}^\intercal) \vct{z}/2} \mathrm{d}\vct{z} \\
&= (2\pi)^{(d-1)^2/2}~Det(\bm{P}\bm{P}^\intercal)^{1/2} Det(\bm{P} \bm{M}^{-1} \bm{P}^\intercal)^{-1/2}. 
\end{align*}

Now it remains to show that $Det(\bm{P} \bm{M}^{-1} \bm{P}^\intercal) = \sum_{T} \prod_{(r,s) \notin T} w_{rs}^{-1} $ where the sum is over all spanning trees of $K_d$. This will finish the proof since we would then have $Det(\bm{P}\bm{P}^\intercal)= \nst(K_d) = 2^{d-1}d^{d-2}$ by Cayley's formula on the number of spanning trees in the complete graph. 

We expand the determinant using the Cauchy-Binet formula. Let $\bm{D} = \bm{M}^{-1/2}$, and let $E$ be the set of edges in $K_d$. For a matrix $\bm{A}$ of size $n \times m$, $I \subseteq \{1,\cdots,n\}, J \subseteq \{1,\cdots,m\}$, we denote by $\bm{A}[I,J]$ the matrix of size $ | I | \times |J|$ whose rows and columns are indexed by $I$ and $J$ respectively. If $I = \{1,\cdots,n\}$, then we write $\bm{A}[~:~,~J]$, and likewise for the column indices. Then, we have
\begin{align}\label{cauchy_binet}
Det(\bm{P}\bm{M}^{-1} \bm{P}^\intercal) = \sum_{\underset{|S| = (d-1)^2}{S \subseteq E}} Det(\bm{P}\bm{D}[~:~,~S~])^2.
\end{align}   
Now we use the following key lemma that we prove later. 
\begin{lemma}\label{invertibility}
Assuming the diagonal entries of the (diagonal) matrix $\bm{D}$ are positive, the matrix $\bm{P}\bm{D}[~:~,~S~]$ is singular if and only if the graph spanned by the complement $\widebar{S} = E \backslash S$ of $S$ in $K_d$ contains a cycle.
\end{lemma}
Since there are exactly $(d-1)$ edges left unchosen by $S$, this lemma implies that they must form a spanning tree in order for the corresponding term to contribute to the sum in identity~\eqref{cauchy_binet}. Hence   
\[Det(\bm{P} \bm{M}^{-1} \bm{P}^\intercal) = \sum_{T~:~\text{spanning tree}} Det(\bm{P}\bm{D}[~:~,~\widebar{T}~])^2.\]
Fix a spanning tree $T$ of $K_d$. Observe that if $T = T^*$ then the edges that generate the cycles in the fundamental cycle basis $\mathcal{B}$ are exactly the ones that are selected in $\widebar{T}$. In other words, each row and each column of $\bm{P}\bm{D}[~:~,~\widebar{T}~]$ contain exactly one non-zero entry, (i.e., $\bm{P}[~:~,~\widebar{T}~]$ is a permutation matrix), hence $Det(\bm{P}\bm{D}[~:~,~\widebar{T}~]) = \pm \prod_{(r,s)\notin T} w_{rs}^{-1/2}$. 
If $T \neq T^*$ then we split the set of edges in $\widebar{T}$ into those that belong to $T^*$ and those that do not. Each column in $\bm{P}\bm{D}[~:~,~\widebar{T}~]$ corresponding to an edge in $\widebar{T} \cap \widebar{T^*}$ contains only one non-zero entry (since this edge is contained in only one cycle in $\mathcal{B}$). Therefore all such edges (columns) along with the corresponding cycles (rows of the non-zero entry) can be successively eliminated from the determinant, yielding
 \begin{equation}
 Det(\bm{P}\bm{D}[~:~,~\widebar{T}~]) = \pm \left(\prod_{(r,s)\in \bar{T} \cap \bar{T^*}} w_{rs}^{-1/2}\right) \cdot Det\left(\bm{P}\bm{D}[~T\cap \widebar{T^*}~,~\widebar{T} \cap T^*~]\right).
 \end{equation}
 Notice that this operation has drastically reduced the size of the problem; the common size $k$ of the sets $T\cap \widebar{T^*}$ and $\widebar{T} \cap T^*$ is anywhere between 0 and $d-1$ at most. Now we will show that 
 \[Det\left(\bm{P}\bm{D}[~T\cap \widebar{T^*}~,~\widebar{T} \cap T^*~]\right) = \pm \prod_{(r,s)\in \bar{T} \cap T^*} w_{rs}^{-1/2},\]
 using a peeling argument slightly more delicate than the one previously applied.
 Observe that $\bm{P}\bm{D}[~T\cap \widebar{T^*}~,~\widebar{T} \cap T^*~]$ is the $k \times k$ cycle-edge incidence matrix with $k$ edges $T\cap \widebar{T^*}$ indexing the rows and $k$ edges in $\widebar{T} \cap T^*$ indexing the columns, such that a row indexed by $e$ indicates the edges $e' \in \widebar{T} \cap T^*$ that participate in the cycle $C_e$.  

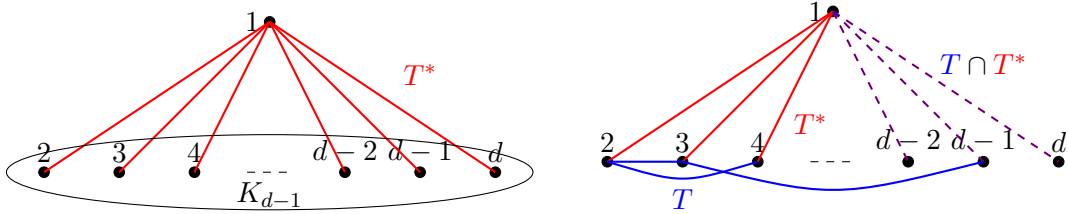
\begin{figure}[h]
\centering
\begin{tikzpicture}

\coordinate [label=left: $1$] (1) at (0,2);

\coordinate [label=above: $2$] (2) at (-3,0);
\coordinate [label=above: $3$] (3) at (-2,0);
\coordinate [label=above: $4$] (4) at (-1,0);

\coordinate [label=above: $d-2$] (5) at (1,0);
\coordinate [label=above: $d-1$] (6) at (2,0);
\coordinate [label=above: $d$] (7) at (3,0);

\filldraw [black] 
(1) circle (2pt) 
(2) circle (2pt)
(3) circle (2pt)
(4) circle (2pt)
(5) circle (2pt) 
(6) circle (2pt)
(7) circle (2pt);
\node (tree) at (2,1.3) {$\textcolor{red}{T^*}$};

\draw[-,thick,red] (1) -- (2); 
\draw[-,thick,red] (1) -- (3);
\draw[-,thick,red] (1) -- (4);
\draw[-,thick,red] (1) -- (5);
\draw[-,thick,red] (1) -- (6);
\draw[-,thick,red] (1) -- (7);
% ...
\draw[-,dashed,thin] (-.3,0) -- (.3,0); 

\node [below,align=center,midway] {$K_{d-1}$};
\draw (0,0) ellipse (3.5cm and .5cm);
 
\end{tikzpicture}
\hspace{.5cm}
\begin{tikzpicture}

\coordinate [label=left: $1$] (1) at (0,2);

\coordinate [label=above: $2$] (2) at (-3,0);
\coordinate [label=above: $3$] (3) at (-2,0);
\coordinate [label=above: $4$] (4) at (-1,0);

\coordinate [label=above: $d-2$] (5) at (1,0);
\coordinate [label=above: $d-1$] (6) at (2,0);
\coordinate [label=above: $d$] (7) at (3,0);

\filldraw [black] 
(1) circle (2pt) 
(2) circle (2pt)
(3) circle (2pt)
(4) circle (2pt)
(5) circle (2pt) 
(6) circle (2pt)
(7) circle (2pt);
\node (tree) at (2,1.3) {$\textcolor{blue}{T} \cap \textcolor{red}{T^*}$};
\node (tree) at (-.3,.5) {$\textcolor{red}{T^*}$};

\draw[-,thick,red] (1) -- (2); 
\draw[-,thick,red] (1) -- (3);
\draw[-,thick,red] (1) -- (4);
\draw[-,thick,dashed,red] (1) -- (5);
\draw[-,thick,dashed,red] (1) -- (6);
\draw[-,thick,dashed,red] (1) -- (7);
% ...
\draw[-,dashed,thin] (-.3,0) -- (.3,0); 

 \node (tree_2) at (-2,-.5) {$\textcolor{blue}{T}$};

\draw[-,thick,blue] (-3,0) -- (-2,0);
\draw[-,thick,blue] (-3,0) .. controls (-2,-.3) .. (-1,0);
\draw[-,thick,blue] (-2,0) .. controls (0,-.5) .. (2,0);

\draw[-,thin,dashed,blue] (1) -- (5);
\draw[-,thin,dashed,blue] (1) -- (6);
\draw[-,thin,dashed,blue] (1) -- (7);

\end{tikzpicture}
\caption{Left: the graph $K_d$ where the star tree $T^*$ is highlighted in red. Right: remaining edges in red and blue after the first elimination procedure (violet edges were removed).}
\label{spanning_tree_1}
\end{figure}

So far, the spanning tree $T^*$ was arbitrary. To continue, we choose $T^*$ to be the \emph{star tree} rooted at vertex 1 (see Figure~\ref{spanning_tree_1}, left). This choice simplifies the combinatorial argument to come, because the fundamental cycle basis $\mathcal{B}$ is now composed of triangles rooted at vertex 1. Crucially, this is where the assumption $G=K_d$ is needed; to ensure the existence of a star spanning tree. Figure~\ref{spanning_tree_1} (right) illustrates the remaining edges after the first elimination procedure.

Since $T$ is a tree, by Lemma 9, each row and column of the matrix $\bm{P}\bm{D}[~T\cap \widebar{T^*}~,~\widebar{T} \cap T^*~]$ contains at least one non-zero entry. Furthermore, $T^*$ being the star graph, each cycle $C_e \in \mathcal{B}$ is a triangle rooted at vertex 1, thus each row of the above matrix contains at most two non-zero entries. This is simply because one of the three edges that compose the triangle $C_e$---namely $e$---is not selected by the set $\widebar{T} \cap T^*$ that indexes the columns of the matrix. See Figure~\ref{spanning_tree_1}, right (any blue edge has at most two adjacent red edges). 

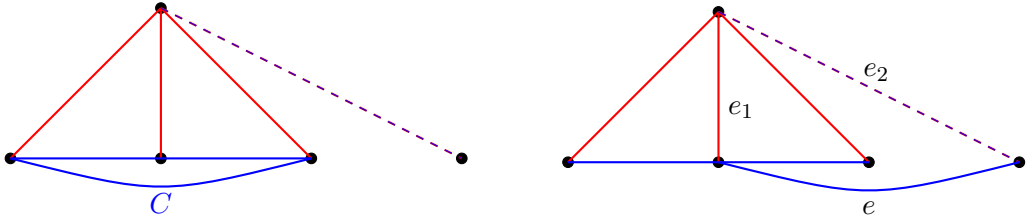
\begin{figure}[h]
\centering
\begin{tikzpicture}

\filldraw [black] 
(0,2) circle (2pt) 

(-2,0) circle (2pt)
(0,0) circle (2pt)
(2,0) circle (2pt)
(4,0) circle (2pt);

\draw[-,thick,red] (0,2) -- (-2,0); 
\draw[-,thick,red] (0,2) -- (0,0);
\draw[-,thick,red] (0,2) -- (2,0);
\draw[-,thick,dashed,red] (0,2) -- (4,0);

\draw[-,thin,dashed,blue] (0,2) -- (4,0);
\draw[-,thick,blue] (-2,0) -- (0,0);
\draw[-,thick,blue] (0,0) -- (2,0);
\draw[-,thick,blue] (-2,0) .. controls (0,-.5) .. (2,0);

\node (cycle) at (0,-.6) {$\textcolor{blue}{C}$};
\end{tikzpicture}
\hspace{1cm}
\begin{tikzpicture}

\filldraw [black] 
(0,2) circle (2pt) 

(-2,0) circle (2pt)
(0,0) circle (2pt)
(2,0) circle (2pt)
(4,0) circle (2pt);

\draw[-,thick,red] (0,2) -- (-2,0); 
\draw[-,thick,red] (0,2) -- (0,0);
\draw[-,thick,red] (0,2) -- (2,0);
\draw[-,thick,dashed,red] (0,2) -- (4,0);

\draw[-,thin,dashed,blue] (0,2) -- (4,0);
\draw[-,thick,blue] (-2,0) -- (0,0);
\draw[-,thick,blue] (0,0) -- (2,0);
\draw[-,thick,blue] (0,0) .. controls (2,-.5) .. (4,0);

\node (edge_e) at (2,-.6) {$e$};
\node (edge_e1) at (.3,.7) {$e_1$};
\node (edge_e2) at (2.1,1.2) {$e_2$};

\end{tikzpicture}
\caption{Left: an impossible situation where there remains a cycle $C$ where no edge was eliminated in the first step. Right: a logical situation where there exist a fundamental cycle $C_e = (e,e_1,e_2)$ with one edge in $T^*$ only, one edge in $T$ only, and one edge in their intersection.}
\label{spanning_tree_2}
\end{figure}

Furthermore, if all the rows contain exactly two non-zero entries then by the pigeonhole principle (since $|T\cap \widebar{T^*}| = |\widebar{T} \cap T^*|$), there will exist three edges in $T\cap \widebar{T^*}$ that form a cycle $C$ (see Figure~\ref{spanning_tree_2}, left). However, we assumed that $T$ is a tree so this cannot happen. Therefore there must exist at least one row in the matrix with exactly one non-zero entry (i.e., there must exist an edge $e \in T\cap \widebar{T^*}$ such that  $C_e = \{e,e_1,e_2\} \in \mathcal{B}$ with $e_1 \in \widebar{T} \cap T^*$ and $e_2 \in T \cap T^*$). See Figure~\ref{spanning_tree_2}. 

\begin{figure}[h]
\centering
\begin{tikzpicture}

\coordinate  (1) at (0,2);

\coordinate (2) at (-3,0);
\coordinate (3) at (-2,0);
\coordinate (4) at (-1,0);

\coordinate (5) at (1,0);
\coordinate  (6) at (2,0);
\coordinate  (7) at (3,0);

\filldraw [black] 
(1) circle (2pt) 
(2) circle (2pt)
(3) circle (2pt)
(4) circle (2pt)
(5) circle (2pt) 
(6) circle (2pt)
(7) circle (2pt);

\draw[-,thick,red] (-3.8,2.2) -- (-3.3,2.2);
\draw[-,thick,blue] (-3.8,1.7) -- (-3.3,1.7);
\node (tree) at (-2.5,2.2) {$\widebar{T} \cap T^*$};
 \node (tree_2) at (-2.5,1.7) {$T \cap \widebar{T^*}$};
\draw[-,thick,red] (1) -- (2); 
\draw[-,thick,red] (1) -- (3);
\draw[-,thick,red] (1) -- (4);
\draw[-,thick,dashed,red] (1) -- (5);
\draw[-,thick,dashed,red] (1) -- (6);
\draw[-,thick,dashed,red] (1) -- (7);
% ...
\draw[-,dashed,thin] (-.3,0) -- (.3,0); 


 \draw[-,thick,blue] (-3,0) -- (-2,0);
\draw[-,thick,blue] (-3,0) .. controls (-2,-.3) .. (-1,0);
\draw[-,thick,blue] (-2,0) .. controls (0,-.5) .. (2,0);

\draw[-,thin,dashed,blue] (1) -- (5);
\draw[-,thin,dashed,blue] (1) -- (6);
\draw[-,thin,dashed,blue] (1) -- (7);

\end{tikzpicture}
\begin{tikzpicture}
\draw[->,thick,black] (-3.8,1) -- (-3.3,1);

\coordinate  (1) at (0,2);
\coordinate (2) at (-3,0);
\coordinate (3) at (-2,0);
\coordinate (4) at (-1,0);

\coordinate (5) at (1,0);
\coordinate  (6) at (2,0);
\coordinate  (7) at (3,0);

\filldraw [black] 
(1) circle (2pt) 
(2) circle (2pt)
(3) circle (2pt)
(4) circle (2pt)
(5) circle (2pt) 
(6) circle (2pt)
(7) circle (2pt);

\draw[-,thick,red] (1) -- (2); 
\draw[-,thick,dashed,red] (1) -- (3);
\draw[-,thick,red] (1) -- (4);
\draw[-,thick,dashed,red] (1) -- (5);
\draw[-,thick,dashed,red] (1) -- (6);
\draw[-,thick,dashed,red] (1) -- (7);
% ...
\draw[-,dashed,thin] (-.3,0) -- (.3,0); 


 \draw[-,thick,blue] (-3,0) -- (-2,0);
\draw[-,thick,blue] (-3,0) .. controls (-2,-.3) .. (-1,0);
\draw[-,thick,dashed,blue] (-2,0) .. controls (0,-.5) .. (2,0);

\draw[-,thin,dashed,blue] (1) -- (5);
\draw[-,thin,dashed,blue] (1) -- (6);
\draw[-,thin,dashed,blue] (1) -- (7);
\end{tikzpicture}
\begin{tikzpicture}
\draw[->,thick,black] (-3.8,1) -- (-3.3,1);

\coordinate  (1) at (0,2);

\coordinate (2) at (-3,0);
\coordinate (3) at (-2,0);
\coordinate (4) at (-1,0);

\coordinate (5) at (1,0);
\coordinate  (6) at (2,0);
\coordinate  (7) at (3,0);

\filldraw [black] 
(1) circle (2pt) 
(2) circle (2pt)
(3) circle (2pt)
(4) circle (2pt)
(5) circle (2pt) 
(6) circle (2pt)
(7) circle (2pt);

\draw[-,thick,dashed,red] (1) -- (2); 
\draw[-,thick,dashed,red] (1) -- (3);
\draw[-,thick,red] (1) -- (4);
\draw[-,thick,dashed,red] (1) -- (5);
\draw[-,thick,dashed,red] (1) -- (6);
\draw[-,thick,dashed,red] (1) -- (7);
% ...
\draw[-,dashed,thin] (-.3,0) -- (.3,0); 


 \draw[-,thick,dashed,blue] (2) -- (3);
\draw[-,thick,blue] (2) .. controls (-2,-.3) .. (4);
\draw[-,thick,dashed,blue] (3) .. controls (0,-.5) .. (6);

\draw[-,thin,dashed,blue] (1) -- (5);
\draw[-,thin,dashed,blue] (1) -- (6);
\draw[-,thin,dashed,blue] (1) -- (7);
\end{tikzpicture}
\begin{tikzpicture}
\draw[->,thick,black] (-3.8,1) -- (-3.3,1);

\coordinate  (1) at (0,2);

\coordinate (2) at (-3,0);
\coordinate (3) at (-2,0);
\coordinate (4) at (-1,0);

\coordinate (5) at (1,0);
\coordinate  (6) at (2,0);
\coordinate  (7) at (3,0);

\filldraw [black] 
(1) circle (2pt) 
(2) circle (2pt)
(3) circle (2pt)
(4) circle (2pt)
(5) circle (2pt) 
(6) circle (2pt)
(7) circle (2pt);

\draw[-,thick,dashed,red] (1) -- (2); 
\draw[-,thick,dashed,red] (1) -- (3);
\draw[-,thick,dashed,red] (1) -- (4);
\draw[-,thick,dashed,red] (1) -- (5);
\draw[-,thick,dashed,red] (1) -- (6);
\draw[-,thick,dashed,red] (1) -- (7);
% ...
\draw[-,dashed,thin] (-.3,0) -- (.3,0); 


 \draw[-,thick,dashed,blue] (2) -- (3);
\draw[-,thick,dashed,blue] (2) .. controls (-2,-.3) .. (4);
\draw[-,thick,dashed,blue] (3) .. controls (0,-.5) .. (6);

\draw[-,thin,dashed,blue] (1) -- (5);
\draw[-,thin,dashed,blue] (1) -- (6);
\draw[-,thin,dashed,blue] (1) -- (7);

\end{tikzpicture}
\caption{An illustration of the peeling process. ``Wedges" with one edge in $T$ only and the other in $T^*$ only are eliminated successively until no edges remain. Violet edges were eliminated in the first step.}
\label{spanning_tree_3}
\end{figure}
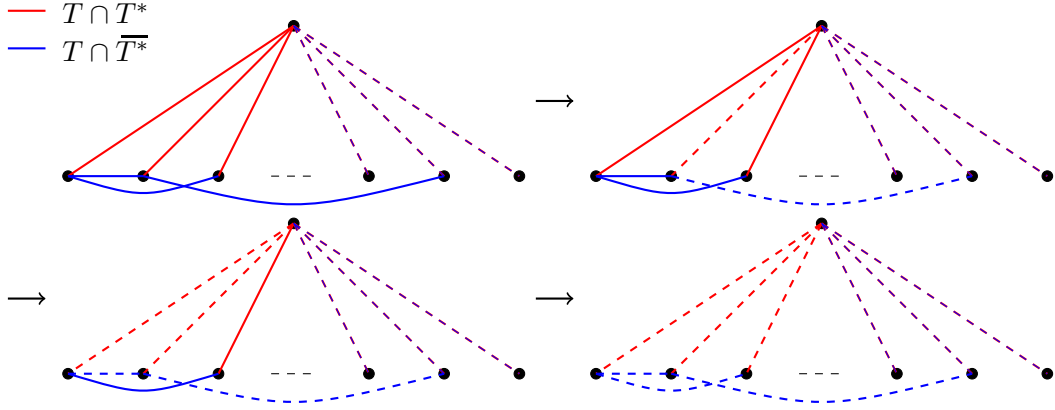

Hence, we can eliminate this row and its corresponding column from the determinant. This corresponds to eliminating (dashing) the edges $e$ and $e_1$ in the right figure above. Applying this argument iteratively allows us to peel all the edges and the cycles they belong to (see Figure~\ref{spanning_tree_3}), so that we obtain   
\[Det\left(\bm{P}\bm{D}[~T\cap \widebar{T^*}~,~\widebar{T} \cap T^*~]\right) = \pm \prod_{(r,s)\in \bar{T} \cap T^*} w_{rs}^{-1/2}.\]
This completes the proof.

\vspace{.5cm}

\begin{proofof}{Lemma~\ref{invertibility}}
Since we assumed the entries of the diagonal matrix $\bm{D}$ are strictly positive, we assume without loss of generality that $\bm{D}$ is the identity matrix. Assume now that the complement of $S$ contains a cycle whose indicator vector is $\vct{\xi} \in \{0,\pm 1\}^{d(d-1)}$. Since $\mathcal{B}$ is a fundamental cycle basis, there exists $\vct{x} \in \R^{(d-1)^2}\backslash \{0\}$ such that $\vct{\xi} = \sum_{e \notin T^*} x_e \vct{\chi}_e = \bm{P}^\intercal \vct{x}$. Since $S$ selects no edges in the cycle indicated by $\vct{\xi}$, it is clear that $\vct{x}^\intercal \bm{P}[~:~,~S~] = 0$, and this settles one direction. As for the other direction, let $\vct{x} \in \R^{(d-1)^2}\backslash \{0\}$ lie in the null space of $(\bm{P}[~:~,~S~])^\intercal$. The vector $\bm{P}^\intercal \vct{x}$ indicates the weights of a Eulerian subgraph in $K_d$ (this vector belong to $\Fspace$ when written in the form of a $d \times d$ matrix). The condition $(\bm{P}[~:~,~S~])^\intercal \vct{x} =0$ implies that this Eulerian subgraph involves no edges from $S$. In particular, any cycle from this subgraph (there always exists one) is in the complement of $E$. This completes the proof.          
\end{proofof}

\subsection{An analytic proof}
This proof contrasts with the previous purely combinatorial approach in that it is mainly analytic. The approach relies on an interpolation argument that involves expressing the Gaussian integral over $\Fspace$ as the \emph{limit} of another parameterized Gaussian integral, when the parameter tends to zero. This latter integral can on the other hand be written in closed form, by relating it to the characteristic polynomial of a Laplacian matrix. Then the Principal Minors Matrix-Tree Theorem is invoked to finish the argument. This final step is the only place where combinatorics appear. (This proof approach was suggested to us by Andrea Sportiello.) Incidentally, this proof can be carried out with an arbitrary graph $G$; there is no need to reduce to the complete case.      
For $\delta >0$ let 
\[I(\delta) = \frac{1}{(2\typeprop \delta^2)^{(d-1)/2}} \int_{\R^{d \times d}} e^{-\frac{1}{2}\sum_{rs} x_{rs}^2/w_{rs}} ~e^{-\frac{1}{2\delta^2}\twonorm{(\x-\x^\intercal)\one}^2} ~\textup{d}\x.\]
The additional Gaussian term in $I(\delta)$ gradually concentrates the mass of the integral on $\Fspace$ as $\delta$ becomes small, and we have the following limiting statement: 
\begin{lemma}\label{limit_integral}
We have 
\[\lim_{\delta \to 0} I(\delta) = c_d ~\int_{\Fspace} ~e^{-\frac{1}{2}\sum_{rs} x_{rs}^2/2w_{rs}} ~\textup{d}\x,\]
with 
\[c_d =  \frac{1}{(2\typeprop)^{(d-1)/2}}\int_{\Fspace^\perp} e^{-2\twonorm{\z\one}^2} ~\textup{d}\z = (2d)^{-(d-1)/2}.\]
\end{lemma}

On the other hand, a straightforward computation allows us to write $I(\delta)$ in closed form:
\begin{lemma}\label{expression_integral} 

Let $G = (V,E)$ be a weighted graph with $V= \{1,\cdots,d\}$, $E= \{(r,s) \in V \times V, r \neq s\}$ where the edges are weighted by the array $\w \in \Rp^{d\times d}$. Let $\bm{L}(\w) \in \R^{d \times d}$ be the Laplacian matrix of $G$. For all $\delta >0$, it holds that
\[ I(\delta) = (2\typeprop)^{((d-1)^2+d)/2} \left(\prod_{r,s} w_{rs}\right)^{1/2} \frac{\delta}{Det\left(\delta^2 \bm{I} +\bm{L}(\w)\right)^{1/2}}.\]
\end{lemma}
Now, by the Principal Minors Matrix-Tree Theorem (see, e.g., \cite{chaiken1982combinatorial}), the characteristic polynomial of the Laplacian matrix of a graph admits the following expansion
\[{Det}\left(x \bm{I} + \bm{L}(\w)\right) = \sum_{F} x^{|\text{roots}(F)|} ~ \prod_{(r,s)\in F}w_{rs},\]
where the sum is over all rooted spanning forests $F$ of the graph.  
We finish the argument by taking a limit in $\delta$: 
\[\delta^{2(d-1)} {Det}\left(\bm{I} + \delta^{-2}\bm{L}(\w)\right) = \delta^{-2} {Det}\left(\delta^2 \bm{I} + \bm{L}(\w)\right) ~  \underset{\delta \to 0}{\xrightarrow{\qquad}} ~ d\sum_{T} \prod_{(r,s) \in T} w_{rs} = (2d)^{d-1} T(\w),\]
since the above limit singles out the rooted spanning forests with exactly one root---i.e., rooted spanning trees---from the characteristic polynomial, and there are $d$ ways of choosing the root of a spanning tree. This exactly leads to the desired identity 
\[\int_{\Fspace} ~e^{-\frac{1}{2}\sum_{rs} x_{rs}^2/2 w_{rs}} ~\textup{d}\x = (2\typeprop)^{((d-1)^2+d)/2} \left(\frac{\prod_{r,s} w_{rs}}{T(\w)}\right)^{1/2}.\]    

\begin{proofof}{Lemma~\ref{limit_integral}}
We decompose $\R^{d \times d}$ into the direct sum $\Fspace \oplus \Fspace^\perp$. It is easy to see that $\Fspace^\perp = \{\z = \vct{\lambda} \one^\intercal - \one \vct{\lambda}^\intercal,~ \vct{\lambda} \in \R^{d}\}$ which is a $(d-1)$-dimensional space. For $\x \in \R^{d \times d}$, let $\y \in \R^{d \times d}$ be its orthogonal projection on $\Fspace$, and $\z = \x-\y$. Therefore $(\x-\x^\intercal)\one = (\z-\z^\intercal)\one = 2\z\one = 2(d\vct{\lambda} - (\one^\intercal \vct{\lambda}) \one)$.
For $\delta >0$, we have
\[I(\delta) = \frac{1}{(2\typeprop \delta^2)^{(d-1)/2}} \int_{F \times F^\perp} e^{-\frac{1}{2}\sum_{r,s} (y_{rs}+z_{rs})^2/w_{rs}} ~e^{-\frac{2}{\delta^2}\twonorm{\z\one}^2} ~\textup{d}\y\textup{d}\z. \] 
We make the change of variables $\z' = \z/\delta$:
\[I(\delta) = \frac{1}{(2\typeprop)^{(d-1)/2}} \int_{\Fspace \times \Fspace^\perp} e^{-\frac{1}{2}\sum_{r,s} (y_{rs}+\delta z'_{rs})^2/w_{rs}} ~e^{-2\twonorm{\z'\one}^2} ~\textup{d}\y\textup{d}\z'. \]
By dominated convergence,
\begin{align*}
\lim_{\delta \to 0} I(\delta) &=  \frac{1}{(2\typeprop)^{(d-1)/2}}  \int_{\Fspace \times \Fspace^\perp} e^{-\frac{1}{2}\sum_{r,s} y_{rs}^2/w_{rs}} ~e^{-2\twonorm{\z\one}^2} ~\textup{d}\y\textup{d}\z\\ 
&=  \frac{1}{(2\typeprop)^{(d-1)/2}} \int_{\Fspace} e^{-\frac{1}{2}\sum_{r,s} y_{rs}^2/w_{rs}} \textup{d}\y ~ \int_{\Fspace^\perp} e^{-2\twonorm{\z\one}^2}\textup{d}\z.
\end{align*}
Moreover,
\[\int_{\Fspace^\perp} e^{-2\twonorm{\z\one}^2}\textup{d}\z = (2d)^{(d-1)/2}\int_{\{\vct{\lambda} \in \R^{d}, \one^\intercal \vct{\lambda} = 0\}} e^{-2d^2\twonorm{\vct{\lambda}}^2} ~\textup{d}\vct{\lambda} = (2\typeprop)^{(d-1)/2} (2d)^{-(d-1)/2},\]
where the pre-factor in the first equality comes from the fact that $\|\z\|^2_{F} = 2d\twonorm{\vct{\lambda}}^2$ for $\z = \vct{\lambda} \one^\intercal - \one \vct{\lambda}^\intercal,~ \vct{\lambda} \in \R^{d},~ \one^\intercal \vct{\lambda} = 0$.  
\end{proofof}

\vspace{.5cm}
\begin{proofof}{Lemma~\ref{expression_integral}}
Let $\delta >0$. We linearize the quadratic term $\twonorm{(\x-\x^\intercal)\one}^2$ in $I(\delta)$ by writing the corresponding Gaussian as the Fourier transform of another Gaussian: $\forall \x \in \R^{d \times d}$,
\[ e^{-\frac{1}{2\delta^2}\twonorm{(\x-\x^\intercal)\one}^2} =  \frac{1}{(2\typeprop)^{d/2}} \int_{\R^d} e^{- \imnb \delta^{-1} \vct{y}^\intercal (\x - \x^{\intercal}) \one - \frac{1}{2}\twonorm{\vct{y}}^2} ~\textup{d}\vct{y},\]
where $\imnb^2 = -1$.
Then
\begin{align*}
I(\delta) = \frac{1}{(2\typeprop \delta^2)^{(d-1)/2}} \frac{1}{(2\typeprop)^{d/2}}  \int_{\R^{d \times d}}\int_{\R^d} e^{-\frac{1}{2}\sum_{rs} x_{rs}^2/w_{rs}} ~ e^{- \imnb \delta^{-1} \vct{y}^\intercal (\x - \x^{\intercal}) \one - \frac{1}{2}\twonorm{\vct{y}}^2} ~\textup{d}\x \textup{d}\vct{y}.
\end{align*} 
We complete the square involving $x_{rs}$ in the exponentiated expression: 
\[-\frac{1}{2}\sum_{r,s} x_{rs}^2/w_{rs} - \imnb \delta^{-1} \vct{y}^\intercal (\x - \x^{\intercal}) \one = 
-\frac{1}{2} \sum_{r,s} \frac{1}{w_{rs}} \left(\left(x_{rs} + \imnb \frac{w_{rs}}{\delta}(y_r - y_s)\right)^2 + \frac{w_{rs}^2}{\delta^2} (y_r - y_s)^2\right).\] 
Then by Fubini's theorem,
\[I(\delta) = \frac{1}{(2\typeprop \delta^2)^{(d-1)/2}} \frac{1}{(2\typeprop)^{d/2}} \int_{\R^d}  e^{-\frac{1}{2}\twonorm{\vct{y}}^2 - \frac{1}{2} \sum_{rs}\frac{w_{rs}}{\delta^2} (y_r - y_s)^2} \int_{\R^{d \times d}} e^{-\frac{1}{2} \sum_{rs} \frac{1}{w_{rs}} \left(x_{rs} + \imnb \frac{w_{rs}}{\delta}(y_r - y_s)\right)^2 } ~\textup{d}\x \textup{d}\vct{y}.\]
The inner integral evaluates to $\left(\prod_{r,s}2\typeprop w_{rs}\right)^{1/2}$. Hence
\begin{align*}
I(\delta) &=  \frac{(2\typeprop)^{(d-1)^2/2}}{\delta^{d-1}} \left(\prod_{r,s} w_{rs}\right)^{1/2}  \int_{\R^d}  e^{-\frac{1}{2}\twonorm{\vct{y}}^2 - \frac{1}{2} \sum_{r,s}\frac{w_{rs}}{\delta^2} (y_r - y_s)^2} ~\textup{d}\vct{y} \\
&= \frac{(2\typeprop)^{((d-1)^2+d)/2}}{\delta^{d-1}} \left(\prod_{r,s} w_{rs}\right)^{1/2} {Det}\left(\bm{I} + \delta^{-2}\bm{L}(\w)\right)^{-1/2},
\end{align*}
where $\bm{L}(\w) \in \R^{d \times d}$ is the Laplacian matrix of the weighted graph $G$.
\end{proofof}  

\section{Discussion} 
Our main result, Theorem~\ref{mainthm}, leaves a gap of essentially a factor of two between $\gamma_{\text{low}}$ and $\gamma_{\text{up}}$. This is a limitation of the methods employed. In particular, the upper bound is likely to be loose due to the lack of concentration of the random variable $\Zpart$ about its mean, and this translates to the possibility of existence of a non-trivial interval inside $[\gamma_{\text{low}},\gamma_{\text{up}}]$ where $\Zpart$ is typically close to 1 while its expectation is exponentially large. A sharper bound could be obtained by computing $\E[|\Zpart - 1|^{1/n}]$, or even, and perhaps less ambitiously, $\E[|\Zpart-1|^\beta]$ for some $0<\beta <1$. The first quantity would correspond to the free energy of the model in the limit; the quantity $|\Zpart - 1|^{1/n}$ is believed to concentrate for large $n$, so taking its logarithm before or after averaging would lead to the same outcome.       

In a different vein, the ``sparse" regime where the sets $S_a$ are of constant size $k$ (exactly or on average) could also be of interest. Here, the relevant scaling is one where $m$ is proportional to $n$. The lower bound argument could be easily extended and yields a bound of $\frac{H(\vct{\pi})}{(d-1)\log k}$. As for the upper bound, one could in principle follow the same first moment strategy, but our analysis breaks in a quite serious fashion, in that none of our asymptotic estimates hold true in this regime.      

\paragraph{Acknowledgments} We benefited from insightful conversations with many people. We thank I-Hsiang Wang for sharing the slides of his talk at ITA. We thank Cris Moore for bringing the work of Achlioptas and Naor~\cite{achlioptas2004two_journal} to our attention. We thank Ross Boczar for his input using Mathematica in the trial-and-error process that led to the discovery of Proposition~\ref{gaussian-integral}. We thank Andrea Sportiello for suggesting the interpolation method that is used in our second proof of Proposition~\ref{gaussian-integral}. We thank Nikhil Srivastava for bringing the survey~\cite{biggs1997algebraic} to our attention. Part of this work was performed during the spring of 2016 when FK and LF were visiting the Simons Institute for the Theory of Computing at UC Berkeley. FK  acknowledges funding from the EU (FP/2007-2013/ERC grant agreement 307087-SPARCS). MJ acknowledges the support of the Mathematical Data Science program of the Office of Naval Research under grant number N00014-15-1-2670.

\bibliographystyle{alpha}
\bibliography{histograms}

\newcommand{\etalchar}[1]{$^{#1}$}
\begin{thebibliography}{BCOH{\etalchar{+}}16}

\bibitem[ACO08]{achlioptas2008algorithmic}
Dimitris Achlioptas and Amin Coja-Oghlan.
\newblock Algorithmic barriers from phase transitions.
\newblock In {\em Foundations of Computer Science, 2008. FOCS'08. IEEE 49th
  Annual IEEE Symposium on}, pages 793--802. IEEE, 2008.

\bibitem[AM04]{achlioptas2004chromatic}
Dimitris Achlioptas and Cristopher Moore.
\newblock The chromatic number of random regular graphs.
\newblock In {\em Approximation, Randomization, and Combinatorial Optimization.
  Algorithms and Techniques}, pages 219--228. Springer, 2004.

\bibitem[AN05]{achlioptas2004two_journal}
Dimitris Achlioptas and Assaf Naor.
\newblock The two possible values of the chromatic number of a random graph.
\newblock {\em Annals of Mathematics}, 162(3):1335--1351, 2005.

\bibitem[BCOH{\etalchar{+}}16]{bapst2016condensation}
Victor Bapst, Amin Coja-Oghlan, Samuel Hetterich, Felicia Ra{\ss}mann, and Dan
  Vilenchik.
\newblock The condensation phase transition in random graph coloring.
\newblock {\em Communications in Mathematical Physics}, 341(2):543--606, 2016.

\bibitem[Big97]{biggs1997algebraic}
Norman Biggs.
\newblock Algebraic potential theory on graphs.
\newblock {\em Bulletin of the London Mathematical Society}, 29(6):641--682,
  1997.

\bibitem[BLM15]{bayati2015universality}
Mohsen Bayati, Marc Lelarge, and Andrea Montanari.
\newblock Universality in polytope phase transitions and message passing
  algorithms.
\newblock {\em Annals of Applied Probability}, 25(2):753--822, 2015.

\bibitem[BMNN16]{banks2016information}
Jess Banks, Cristopher Moore, Joe Neeman, and Praneeth Netrapalli.
\newblock Information-theoretic thresholds for community detection in sparse
  networks.
\newblock In {\em 29th Annual Conference on Learning Theory}, volume~49, pages
  383--416, 2016.

\bibitem[BV04]{boyd02}
Stephen Boyd and Lieven Vandenberghe.
\newblock {\em Convex Optimization}.
\newblock Cambridge University Press, Cambridge, UK, 2004.

\bibitem[Cha82]{chaiken1982combinatorial}
Seth Chaiken.
\newblock A combinatorial proof of the all minors matrix tree theorem.
\newblock {\em SIAM Journal on Algebraic Discrete Methods}, 3(3):319--329,
  1982.

\bibitem[CO09]{coja2009random}
Amin Coja-Oghlan.
\newblock Random constraint satisfaction problems.
\newblock {\em arXiv preprint arXiv:0911.2322}, 2009.

\bibitem[COEH16]{coja2016chromatic}
Amin Coja-Oghlan, Charilaos Efthymiou, and Samuel Hetterich.
\newblock On the chromatic number of random regular graphs.
\newblock {\em Journal of Combinatorial Theory, Series B}, 116:367--439, 2016.

\bibitem[COF14]{coja2014analyzing}
Amin Coja-Oghlan and Alan Frieze.
\newblock Analyzing walksat on random formulas.
\newblock {\em SIAM Journal on Computing}, 43(4):1456--1485, 2014.

\bibitem[COHH16]{coja2016walksat}
Amin Coja-Oghlan, Amir Haqshenas, and Samuel Hetterich.
\newblock Walksat stalls well below the satisfiability threshold.
\newblock {\em arXiv preprint arXiv:1608.00346}, 2016.

\bibitem[COMV09]{coja2009spectral}
Amin Coja-Oghlan, Elchanan Mossel, and Dan Vilenchik.
\newblock A spectral approach to analysing belief propagation for 3-colouring.
\newblock {\em Combinatorics, Probability and Computing}, 18(6):881--912, 2009.

\bibitem[COP16]{coja2016belief}
Amin Coja-Oghlan and Will Perkins.
\newblock Belief propagation on replica symmetric random factor graph models.
\newblock {\em arXiv preprint arXiv:1603.08191}, 2016.

\bibitem[DB70]{debruijn1970asymptotic}
Nicolaas~Govert De~Bruijn.
\newblock {\em Asymptotic Methods in Analysis}.
\newblock Dover Publications, 1970.

\bibitem[DH06]{du2006pooling}
Dingzhu Du and Frank Hwang.
\newblock {\em Pooling Designs and Nonadaptive Group Testing: Important Tools
  for {DNA} Sequencing}, volume~18.
\newblock World Scientific Publishing Company, 2006.

\bibitem[DJM13]{donoho2013information}
David~L Donoho, Adel Javanmard, and Andrea Montanari.
\newblock Information-theoretically optimal compressed sensing via spatial
  coupling and approximate message passing.
\newblock {\em IEEE Transactions on Information Theory}, 59(11):7434--7464,
  2013.

\bibitem[DMO12]{dani2012tight}
Varsha Dani, Cristopher Moore, and Anna Olson.
\newblock Tight bounds on the threshold for permuted k-colorability.
\newblock In {\em Approximation, Randomization, and Combinatorial Optimization.
  Algorithms and Techniques}, pages 505--516. Springer, 2012.

\bibitem[DSS15]{ding2015proof}
Jian Ding, Allan Sly, and Nike Sun.
\newblock Proof of the satisfiability conjecture for large k.
\newblock In {\em Proceedings of the Forty-Seventh Annual ACM on Symposium on
  Theory of Computing}, pages 59--68. ACM, 2015.

\bibitem[DSS16]{ding2016satisfiability}
Jian Ding, Allan Sly, and Nike Sun.
\newblock Satisfiability threshold for random regular nae-sat.
\newblock {\em Communications in Mathematical Physics}, 341(2):435--489, 2016.

\bibitem[FPV15]{feldman2015complexity}
Vitaly Feldman, Will Perkins, and Santosh Vempala.
\newblock On the complexity of random satisfiability problems with planted
  solutions.
\newblock In {\em Proceedings of the Forty-Seventh Annual ACM on Symposium on
  Theory of Computing}, pages 77--86. ACM, 2015.

\bibitem[HLB{\etalchar{+}}01]{Heo1163}
M.~Heo, R.~L. Leibel, B.~B. Boyer, W.~K. Chung, M.~Koulu, M.~K. Karvonen,
  U.~Pesonen, A.~Rissanen, M.~Laakso, M.~I.~J. Uusitupa, Y.~Chagnon,
  C.~Bouchard, P.~A. Donohoue, T.~L. Burns, A.~R. Shuldiner, K.~Silver, R.~E.
  Andersen, O.~Pedersen, S.~Echwald, T.~I.~A. S{\o}rensen, P.~Behn, M.~A.
  Permutt, K.~B. Jacobs, R.~C. Elston, D.~J. Hoffman, and D.~B. Allison.
\newblock Pooling analysis of genetic data: The association of leptin receptor
  ({LEPR}) polymorphisms with variables related to human adiposity.
\newblock {\em Genetics}, 159(3):1163--1178, 2001.

\bibitem[KMZ12]{krzakala2012reweighted}
Florent Krzakala, Marc M{\'e}zard, and Lenka Zdeborov{\'a}.
\newblock Reweighted belief propagation and quiet planting for random {K-SAT}.
\newblock {\em arXiv preprint arXiv:1203.5521}, 2012.

\bibitem[KZ09]{krzakala2009hiding}
Florent Krzakala and Lenka Zdeborov{\'a}.
\newblock Hiding quiet solutions in random constraint satisfaction problems.
\newblock {\em Physical Review Letters}, 102(23):238701, 2009.

\bibitem[MT11]{mezard2011group}
Marc M{\'e}zard and Cristina Toninelli.
\newblock Group testing with random pools: Optimal two-stage algorithms.
\newblock {\em IEEE Transactions on Information Theory}, 57(3):1736--1745,
  2011.

\bibitem[Roc70]{rockafellar70}
R.~Tyrrell Rockafellar.
\newblock {\em Convex Analysis}.
\newblock Princeton University Press, Princeton, 1970.

\bibitem[SBC{\etalchar{+}}02]{sham2002dna}
Pak Sham, Joel~S Bader, Ian Craig, Michael O'Donovan, and Michael Owen.
\newblock {DNA} pooling: a tool for large-scale association studies.
\newblock {\em Nature Reviews Genetics}, 3(11):862--871, 2002.

\bibitem[SSZ16]{sly2016number}
Allan Sly, Nike Sun, and Yumeng Zhang.
\newblock The number of solutions for random regular nae-sat.
\newblock {\em arXiv preprint arXiv:1604.08546}, 2016.

\bibitem[Tan02]{tanaka2002statistical}
Toshiyuki Tanaka.
\newblock A statistical-mechanics approach to large-system analysis of {CDMA}
  multiuser detectors.
\newblock {\em IEEE Transactions on Information theory}, 48(11):2888--2910,
  2002.

\bibitem[Vaa79]{vaaler1979}
Jeffrey~D. Vaaler.
\newblock A geometric inequality with applications to linear forms.
\newblock {\em Pacific Journal of Mathematics}, 83(2):543--553, 1979.

\bibitem[WHLC16]{wang2016data}
I-Hsiang Wang, Shao-Lun Huang, Kuan-Yun Lee, and Kwang-Cheng Chen.
\newblock Data extraction via histogram and arithmetic mean queries:
  Fundamental limits and algorithms.
\newblock In {\em 2016 IEEE International Symposium on Information Theory
  (ISIT)}, pages 1386--1390. IEEE, 2016.

\bibitem[WV09]{wu2009fundamental}
Yihong Wu and Sergio Verd{\'u}.
\newblock Fundamental limits of almost lossless analog compression.
\newblock In {\em 2009 IEEE International Symposium on Information Theory},
  pages 359--363. IEEE, 2009.

\bibitem[Zig04]{zigangirov2004theory}
Kamil~Sh. Zigangirov.
\newblock {\em Theory of Code Division Multiple Access Communication},
  volume~6.
\newblock John Wiley \& Sons, 2004.

\bibitem[ZK15]{zdeborova2015statistical}
Lenka Zdeborov{\'a} and Florent Krzakala.
\newblock Statistical physics of inference: Thresholds and algorithms.
\newblock {\em arXiv preprint arXiv:1511.02476}, 2015.

\bibitem[ZKMZ13]{zhang2013non}
Pan Zhang, Florent Krzakala, Marc M{\'e}zard, and Lenka Zdeborov{\'a}.
\newblock Non-adaptive pooling strategies for detection of rare faulty items.
\newblock In {\em 2013 IEEE International Conference on Communications
  Workshops (ICC)}, pages 1409--1414. IEEE, 2013.

\end{thebibliography}

\end{document}